\documentclass[12pt,a4paper]{amsart}
\usepackage[pagewise]{lineno}
\usepackage[top=35mm, bottom=32mm, left=30mm, right=30mm]{geometry}
\usepackage{newtxtext,newtxmath}
\usepackage{mathrsfs}
\usepackage{fancyhdr} 
\usepackage{url}
\usepackage[all]{xy}
\usepackage[colorlinks=true,citecolor=blue]{hyperref}

\theoremstyle{plain}
\newtheorem{claim}{Claim}
\newtheorem{thm}{Theorem}[section]
\newtheorem{lem}[thm]{Lemma}
\newtheorem{prop}[thm]{Proposition}

\theoremstyle{definition}
\newtheorem{defn}[thm]{Definition}
\newtheorem{rem}[thm]{Remark}
\newtheorem{ques}[thm]{Question}

\newcommand{\Z}{{\mathbb{Z}}}
\newcommand{\N}{\mathbb{N}}

\newcommand{\ep}{\varepsilon}

\DeclareMathOperator{\diam}{diam}
\DeclareMathOperator{\supp}{supp}

\DeclareMathOperator{\sep}{sep}

\newcommand{\pr} {\mathrm{Prob}}

\newcommand{\Map}{\operatorname{Map}}
\newcommand{\Sep}{\operatorname{Sep}}

\newcommand{\vol}{\operatorname{vol}}

\def \d {\delta}

\def \ov {\overline}

\begin{document}

	\title{Relative entropy, topological pressure and variational principle for locally compact sofic group actions}
	
	\author{xianqiang Li and Zhuowei Liu *}
	\date{\today}
	\address[X. Li]{School of Mathematics (Zhuhai), Sun Yat-sen University,
		Zhuhai, Guangdong, 519000, P.R. China}
	\email{lixq233@mail2.sysu.edu.cn}
	\address[Z. Liu]{School of Mathematics (Zhuhai), Sun Yat-sen University,
		Zhuhai, Guangdong, 519000, P.R. China}
	\email{liuzhw55@mail2.sysu.edu.cn}

	\thanks{* Zhuowei Liu is the corresponding author.}
	
	\begin{abstract}
		For a locally compact sofic group continuously acting on a compact metric space, we first study the relative  sofic entropy and prove an additive inequality relating sofic entropy and relative  sofic entropy. Moreover, it is shown that the relative variational principle remains valid in this paper. Secondly, the topological pressure  for locally compact sofic group actions is investigated and the
		variational principle for topological pressure in this soﬁc context is established. As an application, we show a sufficient condition for a signed measure to be a $G$-invariant measure. These contributions generalize the classical results  for countable soﬁc groups on such spaces.
	\end{abstract}
	
	\maketitle
	\section{Introduction}
	
	By a \textit{$G$-topological dynamical system} ($G$-system for short), we mean a pair $(X,G)$, where $X$ is a compact metric space
	and $G$ acts on $X$ as a group of homeomorphsims. For a $G$-system $(X, G)$ if there exists a $G$-invariant Borel probability measure $\mu$ on $X$, it induces a measurable dynamical $G$-system $(X, \mathcal{B}_X, \mu, G)$ , where $\mathcal{B}_X$ is the Borel $\sigma$-algebra on $X$.
	
	Entropy and topological pressure are central concepts in the theory of dynamical systems, providing powerful tools for characterizing complexity, chaotic behavior, and statistical properties. In \cite{R73}, Ruelle introduced topological pressure for $\mathbb{Z}^d$-actions dynamical system  and established the variational principle under expansiveness and specification condition. Walters \cite{W751} later proved the principle for $\mathbb{Z}_+$-actions without these assumptions, and a shorter proof for $\mathbb{Z}_+^d$-actions was given by Misiurewicz \cite{M76}. The classical variational principle establishes a profound connection between topological pressure and measure-theoretic entropy, stating that the topological pressure is the supremum of the sum of the measure-theoretic entropy and the free energy over all invariant measures, and that this supremum is attained under certain conditions. 
	
	The classical theory of entropy and topological pressure  primarily studies actions generated by the integer group  $\Z$ or its powers $\Z^d$. A significant direction in the modern development of dynamical system and ergodic theory is the extension of this framework from classical  $\Z$ or  $\Z^d$ actions to actions of more general group, such as countable discrete group, compact group, and locally compact group. This generalization faces fundamental challenges, as many desirable properties that hold for \(\mathbb{Z}\)-actions—such as the existence of invariant measures and the feasibility of averaging processes—may fail for actions of more general group. The entropy theory for actions of amenable group is relatively well-developed, benefiting from the presence of Følner sequences, which provide natural finite approximations for averaging(see for example \cite{D01,K75,OW87,RW00}). The variational principle for topological pressure was  extended to actions of any countable amenable group by Ollagnier and Pinchon \cite{OP82,O85}, Stepin and Tagi-Zade \cite{ST80}, and Tempelman \cite{T84,T92}.
	
	However, for non-amenable groups (e.g., free groups), the classical averaging-based definition of entropy is no longer applicable, prompting the development of new theories such as sofic entropy. The notion of sofic groups was first introduced by Gromov \cite{GM} and explicitly by Weiss in \cite{W00}, which form a broad class that includes all amenable groups and many non-amenable groups; their definition relies on the existence of a sequence of permutation groups that asymptotically approximate the group action (a sofic approximation). In \cite{Bow10b}, Bowen pioneered the study of entropy for measure-preserving actions of countable sofic group on probability spaces. Extending this framework, Kerr and Li \cite{KL11b,KERRLIHD} introduced the notions of topological entropy  and further established a variational principle bridging the measure-theoretic and topological perspectives. This achievement extended the variational principle to a very general framework, marking a milestone in the entropy theory for actions of countable group. In \cite{C13}, N. Chung generalized it to the variational principle of topological pressure under countable sofic group actions.

	On the other hand, relative entropy (or conditional entropy) represent another crucial refinement in the classical theory. Let $(X,G)$ and $(Y,G)$ be two $G$-systems. A continuous map $\pi: X \to Y$ is called  a \textit{factor map} between $(X,G)$ and $(Y,G)$ if it is onto and $\pi(gx) = g\pi(x)$ for all $g \in G$ and $x \in X$. In this case, we say $(X,G)$ is an extension of $(Y,G)$ or $(Y,G)$ is a factor of $(X,G)$. The relative entropy aim to describe the complexity of a system relative to a given  factor system. This may more deeply analyze the complexity of the fibers or the factor's $\sigma$-algebra.  The authors in \cite{DS02,LW93} proved the relative version of variational principle. Recent work by \cite{Y15,ZZC15} introduced the notion of conditional amenable entropy and established a corresponding variational principle. This was extended by Luo \cite{luo17}, who studied the conditional sofic entropy and relative sofic entropy of $Y$ with respect to $X$ for both topological and measurable dynamical systems, generalizing the relative variational principle to countable sofic groups setting. In \cite{G21}, Hayes futher investigated the relative sofic entropy and proved the outer Pinsker factor of a product action is the product of the outer Pinsker factors in many cases.

	Despite the remarkable success of sofic theory in dealing with discrete non-amenable groups, its focus has primarily been on actions of discrete group. However, in many problems from mathematical physics and geometry (e.g., actions of Lie group on homogeneous spaces, continuous-time stochastic processes), actions of locally compact group (such as \(\mathbb{R}^n\), \(\operatorname{SL}(n, \mathbb{R})\)) are more natural and prevalent. Extending entropy theory from discrete to locally compact groups presents essential new challenges: the structure of a locally compact group is no longer discrete, and its "size" is characterized by a Haar measure. A key obstacle is how to integrate discrete sofic approximations with continuous group actions while properly handling analytical issues such as compactness and measure continuity. In \cite{BB22}, Bowen and Burton presented a novel approach to locally compact sofic groups, based on partial actions and charts. Subsequently, in \cite{B25}, Bowen studied the locally compact sofic entropy theory and also established a variational principle similar to the discrete case. (This theoretical framework was initially introduced by Sukhpreet Singh in his unpublished PhD thesis \cite{S16} under Bowen's direction.)
	
	Thus it is a very natural question whether the relative entropy, topological pressure and variational principle can be generalized to locally compact sofic group actions. This  paper tries to use the method of defining locally compact sofic groups described in \cite{B25,BB22} to address this question. To be precise, we have defined the relative entropy in two different ways and proved their equivalence(see Proposition~\ref{eq}). This allows us to establish an additive inequality relating entropy and relative entropy(see Theorem~\ref{add}). Furthermore, we demonstrate a relative version of the variational principle for topological entropy and measure entropy(see Theorem~\ref{vxy}). On the other hand, we introduce the topological pressure for actions of locally compact sofic group and prove a variational principle for it(see Theorem~\ref{main3}). Additionally, we investigate fundamental properties of topological pressure. As an application, we provide a sufficient condition for a signed measure to be a $G$-invariant measure(see Theorem~\ref{signed}).
	
	The paper is organized as follows. In Section \ref{sec2}, we recall some definitions and some related theorems. In 
	Section \ref{sec3}, we study the relative sofic topological entropy and prove the additive inequality.  In Sections \ref{sec4}, we investigate the relative sofic measure entropy and show the relative version of the variational principle. In Section \ref{sec5}, we introduce the notion of topological pressure for actions of locally compact sofic group. Then we prove the variational principle for topological pressure and study the fundamental properties of topological pressure.
	
	\textit{Acknowledgments}.
	The authors would like to thank their doctoral supervisor, Siming Tu, who provided much guidance and assistance.
	
	\section{Preliminaries}\label{sec2}
	\subsection{Topological dynamical systems and factor maps}
	Let $G$ be a locally compact second countable group with the identity $1_G$. By
	a \emph{$G$-topological dynamical system} ($G$-system for short) $(X, G)$
	we mean that $X$ is a compact metric space and $\Pi \colon G \times X \to X$, given by $(g, x) \mapsto gx$, is a continuous mapping satisfying $\Pi(1_G, x) = x$ for every $x \in X$ and $\Pi(g_1, \Pi(g_2, x)) = \Pi(g_1 g_2, x)$ for every $g_1, g_2 \in G$ and $x \in X$.	Let $(X,G)$ and $(Y,G)$ be two $G$-systems. A continuous map $\pi: X \to Y$ is called a \textit{homomorphism} or a \textit{factor map} between $(X,G)$ and $(Y,G)$ if it is onto and $\pi(gx) = g\pi(x)$ for all $g \in G$ and $x \in X$. In this case we say $(X,G)$ is an \textit{extension} of $(Y,G)$ or $(Y,G)$ is a \textit{factor} of $(X,G)$. If $\pi$ is also injective then it is called an \textit{isomorphism}.

	\subsection{Local $G$-spaces and locally compact sofic groups}
	
	\begin{defn}\label{defn}
		A \emph{partial left-action} of $G$ on a Hausdorff space $M$ is a continuous map $\alpha \colon \operatorname{dom}(\alpha) \to M$ where $\operatorname{dom}(\alpha) \subset G \times M$ is open. We require the following axioms hold for all $p \in M$:
		\begin{enumerate}
			\item \textbf{Axiom 1:} $(1_G, p) \in \operatorname{dom}(\alpha)$ and $\alpha(1_G, p) = p$.
			\item \textbf{Axiom 2:} If $(g, p) \in \operatorname{dom}(\alpha)$ then $(g^{-1}, \alpha(g, p)) \in \operatorname{dom}(\alpha)$ and $\alpha(g^{-1}, \alpha(g, p)) = p$.
			\item \textbf{Axiom 3:} If $(h, p), (g, \alpha(h, p)), (gh, p) \in \operatorname{dom}(\alpha)$ then $\alpha(gh, p) = \alpha(g, \alpha(h, p))$.
		\end{enumerate}
		A partial action $\alpha$ is \emph{homogeneous} if it additionally satisfies:
		\begin{enumerate}
			\setcounter{enumi}{3}
			\item \textbf{Axiom 4:} For every $p \in M$ there is an open neighborhood $O_p$ of $1_G$ in $G$ such that $O_p \times \{p\} \subset \operatorname{dom}(\alpha)$ and the restriction of $\alpha(\cdot, p)$ to $O_p \times \{p\}$ is a homeomorphism onto an open neighborhood of $p$ in $M$.
		\end{enumerate}
	\end{defn}
	A \emph{local left-$G$-space} is a pair $(M, \alpha)$ where $M$ is an Hausdorff space and $\alpha$ is a partial homogeneous left-action of $G$.

	\noindent\textbf{Notation.} Throughout this paper, a pair $(M, \alpha)$ as defined above is called a local $G$-space, and will typically be denoted simply by $M$, with the action $\alpha$ understood. We adopt the concise notation $g.p$ for the partial left action $\alpha(g,p)$. For a subset $K \subset M$, we write $g.K = \{ g.k :k \in K \}$, which is defined precisely when $\{g\} \times K \subset \operatorname{dom}(\alpha)$. Likewise, for a subset $O \subset G$ and a point $p \in M$, the orbit $O.p = \{ g.p : g \in O \}$ is defined when $O \times \{p\} \subset \operatorname{dom}(\alpha)$.

	\begin{defn}
		Let $(M, \alpha)$ be a local $G$-space and take a point $p \in M$. A \emph{chart centered at $p$} is a homeomorphism $f_p : \operatorname{dom}(f_p) \to \operatorname{rng}(f_p)$ satisfying the following conditions:
		\begin{itemize}
			\item The domain $\operatorname{dom}(f_p)$ is an open neighborhood of $p$ in $M$;
			\item The range $\operatorname{rng}(f_p)$ is an open neighborhood of the identity element in the group $G$;
			\item For every $g \in \operatorname{rng}(f_p)$, we have $g = f_p(g. p)$, meaning the chart $f_p$ maps the result of the partial left action back to the group element itself.
		\end{itemize}
		We require that the partial left action $g.p$ is well-defined for all $g \in \operatorname{rng}(f_p)$. According to Axiom 4 of Definition \ref{defn}, for every point $p \in M$, there exists such a chart centered at $p$.
	\end{defn}
	
	The following two results are from \cite{BB22}, which show that there is a canonical measure with locally measure-preserving property  on local $G$-space.
	
	\begin{prop}
		Let $(M, \alpha)$ be a local $G$-space. Fix a right-Haar measure $\mathrm{Haar}$ on $G$. Then there exists a unique Radon measure $\mathrm{vol}_M$ on $M$ satisfying the following. If $p \in M$, $f_p$ is a chart centered at $p$ and $K \subset \operatorname{dom}(f_p)$ is Borel then
		\begin{align*}
			\mathrm{vol}_M(K) = \mathrm{Haar}(\{g \in \operatorname{rng}(f_p) : g.p \in K\}) = \mathrm{Haar}(f_p(K)).
		\end{align*}
		We write $\mathrm{vol}$ instead of $\mathrm{vol}_M$, when the choice of $M$ is clear.
	\end{prop}
	
	\begin{lem}\label{lmp}
		Let $(M, \alpha)$ be a local $G$-space and suppose $\{g\} \times K \subset \operatorname{dom}(\alpha)$ for some measurable $K \subset M$ and $g \in G$. If $G$ is unimodular then
		\[
		\mathrm{vol}(g.K) = \mathrm{vol}(K).
		\]
	\end{lem}
	By \cite{BB22}, any soﬁc group is unimodular. 
	
	\begin{defn}
		Let $\mathcal{M} = (M, \alpha)$ be a local $G$-space, $U$ an open precompact subset of $G$, and $\epsilon > 0$. Define $M[U] = M[\alpha, U]$ as the set of all $p \in M$ such that:
		\begin{enumerate}
			\item If $g, h \in G$ satisfy $g, h, gh \in U$, then $g .(h .p) = (gh) . p$ (with both sides well-defined)
			\item The map $g \mapsto \alpha(g, p)$ is a homeomorphism from $U$ to an open neighborhood of $p$
		\end{enumerate}
		We say $\mathcal{M}$ is a \emph{$(U, \epsilon)$-sofic approximation} to $G$ if $\operatorname{vol}(M) < \infty$ and
		\[
		\operatorname{vol}(M[U]) \geq (1 - \epsilon) \operatorname{vol}(M).
		\]
	\end{defn}
	
	\begin{defn}
		A \emph{sofic approximation} to $G$ is a sequence $\Sigma = (M_i)_{i=1}^{\infty}$ where each $M_i$ is a $(U_i, \epsilon_i)$-sofic approximation such that:
		\begin{enumerate}
			\item The $U_i$ are precompact open sets increasing to $G$ ($U_1 \subset U_2 \subset \cdots$ and $\bigcup_i U_i = G$);
			\item $\epsilon_i \to 0$ as $i \to \infty$.
		\end{enumerate}
		We say $G$ is \emph{sofic} if it admits a sofic approximation.
	\end{defn}
	In this paper, we always assume $G$ is a locally compact sofic group, $X$ is a compact metric space  and $(X, G)$ always denotes a $G$-system unless otherwise specified.

	\subsection*{3.3. Pseudometrics on model spaces}

	A \emph{pseudometric} $\rho \colon X \times X \to [0, \infty)$ on a set $X$ is a metric on $X$ except that $\rho(x, y) = 0$ does not necessarily imply $x = y$. Let $(X,G)$ be a $G$-system. A  continuous pseudometric $\rho$ on $X$ is called \emph{dynamically generating} if for every $x \neq y$ there exists $g \in G$ such that $\rho(gx, gy) > 0$; is called \emph{1-bounded} if $\rho(x, y) \leq 1$ for all $x, y \in X$.

	\begin{defn}
		For a pseudometric space $(X, \rho)$, a subset $Y \subset X$ is \emph{$(\rho, \epsilon)$-separated} if $\rho(x, y) > \epsilon$ for all distinct $x, y \in Y$; define $\Sep_\epsilon(X, \rho) = \max \{ |Y| : Y \subset X \text{ is } (\rho, \epsilon)\text{-separated} \}$.
		
	\end{defn}

	Given a finite measure space $(M, \vol)$ we let $\overline{\vol}$ denote the probability
	measure $\overline{\vol}(E) = 
	\vol(E)/\vol(M)$ .

	\begin{defn}
		Suppose $(M, \vol)$ is a finite measure space and $X$ is a Hausdorff space. Let $\Map(M, X)$ denote the space of all measurable maps from $M$ to $X$. If $\rho$ is a pseudometric on $X$, define the pseudometric $\rho^M$ on $\Map(M, X)$ by
		\[
		\rho^M(\varphi, \psi) = \vol(M)^{-1} \int_M \rho(\varphi(p), \psi(p))  d\vol(p) = \int_M \rho(\varphi(p), \psi(p))  d\overline{\vol}(p).
		\]
		
		More generally, for measurable maps $\varphi, \psi \colon M \to X$ whose domains may not be all of $M$, let $D \subset M$ be the intersection of their domains. Then define
		\[
		\rho^M(\varphi, \psi) = \overline{\vol}(M \setminus D) + \int_{p \in D} \rho(\varphi(p), \psi(p))  d\overline{\vol}(p).
		\]
	\end{defn}

	\begin{defn}
		Let $M$ be a finite volume local $G$-space and $\varphi \colon M \to X$. For $g \in G$ define $\varphi \circ g$ by 
		\[
		(\varphi \circ g)(p) = \varphi(g . p)
		\] 
		for every $p \in M$ such that $g . p$ is well-defined. We also define $g \circ \varphi \colon M \to X$ associatively by
		\[
		(g \circ \varphi)(p) = g(\varphi(p)).
		\]
	\end{defn}
	
	\begin{defn}
		Let $U \subset G$ be a precompact neighborhood of identity, $\delta > 0$. Suppose $M$ is a finite volume local $G$-space. A map $\varphi \colon M \to X$ is \emph{$(U, \delta, \rho^M)$-equivariant} if
		\[
		\rho^M(\varphi \circ g,\  g \circ \varphi) < \delta
		\]
		for all $g \in U$. Denote by $\Map(M, X, \rho \colon U, \delta) \subset \Map(M, X)$ the set of all $(U, \delta, \rho_M)$-equivariant maps.
	\end{defn}
	

	\section{Relative sofic topological entropy}\label{sec3}
	In this section, we introduce the notion of the relative sofic topological entropy in two different ways and prove their equivalence. Then we show an additive inequality relating sofic topological entropy and relative sofic topological entropy.
	
	\begin{defn}\label{def1}
		Let $\pi:(X,G)\to (Y,G)$ be a factor map between $G$-systems. Let \(\rho_X\) and \(\rho_Y\) be continuous pseudometrics on \(X\) and \(Y\) respectively, $\Sigma = (M_i)_{i=1}^{\infty}$ a sofic approximation to $G$.  For \(\Phi \subset \Map(M, X)\), denote by \(N_\epsilon(\Phi; \rho_{X}^M | Y)\) the maximal cardinality of a \((\rho_{X}^M, \epsilon)\)-separated subset \(\Phi_0 \subset \Phi\) satisfying $\pi \circ \varphi=\pi \circ \psi$,
		for all \(\varphi, \psi \in \Phi_0\). Define
		\begin{align*}
			h_{\Sigma}^\ep(\rho_X|Y) &= \inf_{U \subset G} \inf_{\delta > 0} \limsup_{i \to \infty} \frac{1}{\vol(M_i)} \log N_\epsilon \big( \Map(M_i, X, \rho_{X} \colon U,\d); \rho_{X}^{M_i} | Y\big), \\
			h_{\Sigma}(\rho_X|Y) &= \sup_{\epsilon > 0} h_{\Sigma}^\ep(\rho_X|Y).
		\end{align*}
		where the first infimum is over all precompact neighborhoods $U$ of the identity in $G$. 
	\end{defn}
	When $(Y,G)$ is trivial, we recover standard sofic topological entropy:
	$$	h_{\Sigma}(\rho_X):=h_{\Sigma}(\rho_X|Y).$$
	
	\begin{lem}\label{usefull}
		Let $(X,G)$ be a $G$-system and $\rho_1$ and $\rho_2$ be two 1-bounded dynamically generating continuous pseudometrics on $X$.	Let $U_1\subset G$ be a precompact neighborhood of the identity and $\ep,\d_1>0$. Then for a sufficiently good sofic approximation $M$, there exists a precompact neighborhood of the identity $U_2\subset G$ and $\d_2>0$ such that $$\Map(M, X, \rho_2 \colon U_2,\d_2)\subset\Map(M, X, \rho_1 \colon U_1,\d_1).$$
	\end{lem}
	\begin{proof}
		Let	$f : G \to (0, \infty)$ be a continuous function with $1 = \int_G f(g) \, \mathrm{dHaar}(g) < \infty$.
		Define $\rho_f$ by
		$\rho_1^f (x, y) := \int_G \rho_1(gx, gy) f(g) \, \mathrm{dHaar}(g)$ and similarly define $\rho_2^f (x, y) := \int_G \rho_2(gx, gy) f(g) \, \mathrm{dHaar}(g)$. It is easy to check that $\rho_1^f (x, y)$ and $\rho_2^f (x, y)$ are 1-bounded pseudometrics. Now we show that $\rho_1^f $ and $\rho_2^f $ are continuous metrics.
		Assume that $x,y\in X$ with $\rho_1^f (x, y)=0$. Since $\int_G \rho_1(gx, gy) f(g) \, \mathrm{dHaar}(g)\ge 0$, $\rho_1(gx, gy) f(g)=0$ for a.e. Then $\rho_1(gx, gy)=0$ for a.e. by $f>0$. Particularly, $\rho_1(gx, gy)=0$ for any $g$ in a dense subset of $G$. Since $\Pi \colon G \times X \to X$, given by $(g, x) \mapsto gx$, is a continuous mapping and $\rho_1$ is continuous, the map $g\mapsto \rho_1(gx, gy)$ is also continuous. Then one has $\rho_1(gx, gy)=0$ for any $g\in G$.
		Since $\rho_1$ is dynamically generating, we have $x=y$. Therefore $\rho_1^f$ is a continuous metric. Similarly,  $\rho_2^f$ is also a continuous metric.
		
		We claim that for any given $\xi>0$, there exists a small enough $\eta$ such that for any $x, y \in X$ if $\rho_1^f(x, y) \leq \eta$ then $\rho_2(x, y) \leq \xi$. Suppose this does not hold true. There exists $\xi_0>0$ such that for any $n\in \N$, there exist $x_n,y_n\in X$, $\rho_1^f(x_n, y_n) \leq \frac{1}{n}$ but $\rho_2(x_n, y_n) \ge \xi_0$. Since $X$ is compact, without loss of generality, we assume $(x_n,y_n)$ converges to $(x,y)$ as $n\to \infty$. By the continuity of $\rho_1$ and $\rho_2^f$, $\rho_1^f(x,y)=\lim_{n\to\infty}\rho_1^f(x_n,y_n)\leq\lim_{n\to\infty}\frac{1}{n}=0$. But $\rho_2(x,y)=\lim_{n\to\infty}\rho_2(x_n,y_n)\ge\xi_0$. A contradiction. This completes the proof of the claim.
		
		Then by the claim, let $0 < \eta \leq \ep/2$ be such that for any $x, y \in X$ if $\rho_1^f(x, y) \leq \eta$ then $\rho_2(x, y) \leq \epsilon/4$; and similarly if $\rho_2^f(x, y) \leq \eta$ then $\rho_1(x, y) \leq \delta_1/2$.
		Choose a compact neighborhood of identity $W \subset G$ with $\int_W f(g) \, \mathrm{dHaar}(g) \geq 1 - \eta^2/2$.
		Let $U_2 \subset G$ be a precompact neighborhood of the identity with $WU_1 \subset U_2$ and set $\delta_2 = \min \left\{ \frac{\delta_1^2 \eta}{8}, \frac{\epsilon \eta}{8} \right\}$.
		Now, if $\varphi \in \Map(M, X, \rho_2 \colon U_2,\d_2)$ and $h \in U_1$, we have

		\begin{align*}
			\int_M \int_W \rho_2( g \circ \varphi(h.p), gh \circ \varphi(p) ) f(g) \, \mathrm{dHaar}(g) \, d\overline{\vol}(p)\\ 
			\leq \int_M \int_W \rho_2( g \circ \varphi(h.p), \varphi (gh.p) ) f(g) \, \mathrm{dHaar}(g) \, d\overline{\vol}(p)\\ 
			+ \int_M \int_W \rho_2( \varphi (gh.p), gh \circ \varphi(p) ) f(g) \, \mathrm{dHaar}(g) \, d\overline{\vol}(p) \\ 
			\leq \int_W (\delta_2 + \delta_2) f(g) \,\mathrm{dHaar}(g) \leq 2\delta_2.
		\end{align*}
		
		Let $A := \left\{ p \in M : \int_W \rho_2( g \circ \varphi(h.p), gh \circ \varphi(p) ) f(g) \, \mathrm{dHaar}(g)\leq \frac{\eta}{2} \right\}.$
		By Markov’s inequality, we  have  $\overline{\vol}(M\setminus A) < \frac{4\delta_2}{\eta}$.
		Moreover, observe from our choice of $W$, for every $p \in A$, $\rho_2^f( \varphi \circ (h.p), h \circ \varphi(p) ) \leq \eta$; and thus $\rho_1( \varphi \circ (h.p), h \circ \varphi(p) ) \leq \delta_1/2$.
		Hence for every $h \in U_1$ we conclude that
		\[
		\rho_1^M( \varphi \circ h, h\circ \varphi ) =  \int_M \rho_1( \varphi \circ (h.p), h \circ \varphi(p) \, d\overline{\vol}(p) \leq \frac{\delta_1}{2} + \frac{4\delta_2}{\eta} \leq \delta_1.
		\]
		Therefore,$$\Map(M, X, \rho_2 \colon U_2,\d_2)\subset\Map(M, X, \rho_1 \colon U_1,\d_1).$$
		
	\end{proof}

	\begin{thm}\label{main1}
		Let $\Sigma = (M_i)_{i=1}^{\infty}$ be a sofic approximation to $G$,  $\rho_1$ and $\rho_2$ be two 1-bounded dynamically generating continuous pseudometrics on $X$. Then for any factor map  $\pi:X\to Y$ between two $G$-systems $(X,G)$  and $(Y,G)$ , we have 
		$$h_{\Sigma}(\rho_1|Y)=h_{\Sigma}(\rho_2|Y).$$
	\end{thm}
	\begin{proof}
		
		Let $U_1\subset G$ be a precompact neighborhood of the identity and $\ep,\d_1>0$. Take the same $f$, $\eta$ and $W$ in Lemma~\ref{usefull}. By Lemma~\ref{usefull},  for a sufficiently good sofic approximation $M$, there exists a precompact neighborhood of the identity $U_2\subset G$ and $\d_2= \min \left\{ \frac{\delta_1^2 \eta}{8}, \frac{\epsilon \eta}{8} \right\}>0$ such that $$\Map(M, X, \rho_2 \colon U_2,\d_2)\subset\Map(M, X, \rho_1 \colon U_1,\d_1).$$
		
		Set $\ep'=\ep\eta/4$. Next we show any $(\rho_2^M,\ep)$-separated set in $\Map(M, X, \rho_2 \colon U_2,\d_2)$ is
		$(\rho_1^M,\ep')$-separated in $\Map(M, X, \rho_1 \colon U_1,\d_1)$. So, suppose $\varphi,\psi\in\Map(M, X, \rho_2 \colon U_2,\d_2)$ are such that
		$\rho^M_1(\varphi,\psi)<\ep'$. Then we have
		\begin{align*}
			\int_M \int_W \rho_1( g \circ \varphi(p), g\circ \psi(p) ) f(g) \, \mathrm{dHaar}(g) \, d\overline{\vol}(p)\\ 
			\leq \int_M \int_{G\setminus W} \rho_1( g \circ \varphi(p), g\circ\psi (p) ) f(g) \, \mathrm{dHaar}(g) \, d\overline{\vol}(p)\\ + \int_M \int_W \rho_1( g \circ \varphi(p), g\circ\psi (p) ) f(g) \, \mathrm{dHaar}(g) \, d\overline{\vol}(p) \\ 
			\leq\frac{\eta^2}{2}+\int_M \int_W \rho_1( g \circ \varphi(p),\varphi (g.p) ) f(g) \, \mathrm{dHaar}(g) \, d\overline{\vol}(p) \\
			+ \int_M \int_W \rho_1( \varphi (g.p), \psi (g.p) (p) ) f(g) \, \mathrm{dHaar}(g) \, d\overline{\vol}(p) \\
			+ \int_M \int_W \rho_1( \psi (g.p), g\circ\psi (p) ) f(g) \, \mathrm{dHaar}(g) \, d\overline{\vol}(p) \\
			\leq \frac{\eta^2}{2}+\delta_2+\rho_1^M(\varphi,\psi)+\delta_2 \leq \frac{3\ep\eta}{4}.
		\end{align*}

		Let	$B := \left\{ p \in M : \rho_1^f (\varphi(p), \psi(p)) = \int_G \rho_1(g\circ\varphi(p), g\circ\psi(p)) f(g) \, \mathrm{dHaar}(g)\leq \eta\right\}.$ By Markov’s inequality, we  have  $\overline{\vol}(M\setminus B) \leq\frac{3\ep}{4}$.
		Moreover, for every $p \in B$, one has $\rho_1^f( \varphi (p),\psi(p) ) \leq \eta$; and thus $\rho_2( \varphi(p), \psi(p) ) \leq \ep/4$.
		Hence  we conclude that
		\[
		\rho_2^M( \varphi, \psi ) =  \int_B \rho_2( \varphi  (p), \psi(p) \, d\overline{\vol}(p)+\int_{M\setminus B} \rho_2( \varphi  (p), \psi(p) \, d\overline{\vol}(p) \leq \frac{\ep}{4} + \frac{3\ep}{4} =\ep.
		\]
		Then any \((\rho^{M}_2, \ep)\)-separated subset \(\Phi_0 \subset	\Map(M, X, \rho_2 \colon U_2,\d_2)\) satisfying $\pi \circ \varphi=\pi \circ \psi$
		for all \(\varphi, \psi \in \Phi_0\), is a \((\rho^{M}_1, \ep')\)-separated subset \(\Phi_0 \subset\Map(M, X, \rho_1 \colon U_1,\d_1)\) satisfying $\pi \circ \varphi=\pi \circ \psi$,
		for all \(\varphi, \psi \in \Phi_0\).
		Thus,
		$$N_{\ep} \big( \Map(M, X, \rho_2 \colon U_2,\d_2); \rho^{M}_2 | Y\big)\leq  N_{\ep'}\big( \Map(M, X, \rho_1 \colon U_1,\d_1); \rho^{M}_1 | Y\big).$$
		\begin{align*}
			\limsup_{i \to \infty} \frac{1}{\vol(M_i)} \log N_{\ep} \big( \Map(M_i, X, \rho_2 \colon U_2,\d_2); \rho^{M_i}_2 | Y\big)\leq \\ \limsup_{i \to \infty} \frac{1}{\vol(M_i)} \log N_{\ep'}  \big( \Map(M_i, X, \rho_1 \colon U_1,\d_1); \rho^{M_i}_1| Y\big).
		\end{align*}
		Since $U,\d_1$ and $\ep$ are arbitrary, we have
		$$h_{\Sigma}(\rho_2|Y)\leq h_{\Sigma}(\rho_1|Y).$$
		By symmetry, 	$h_{\Sigma}(\rho_2|Y)\ge h_{\Sigma}(\rho_1|Y).$ Then $h_{\Sigma}(\rho_2|Y)=h_{\Sigma}(\rho_1|Y).$
		
	\end{proof}
	
	We shall denote this common value by $h_\Sigma(X|Y)$, and call it the \emph{relative sofic topological entropy} of $(X,G)$ with respect to the factor $(Y,G)$.
	
	We can also define relative topological entropy by the following method.
	
	\begin{defn}
		Let $\pi:(X,G)\to (Y,G)$ be a factor map between $G$-systems. Let \(\rho_X\) and \(\rho_Y\) be continuous pseudometrics on \(X\) and \(Y\) respectively, $\Sigma = (M_i)_{i=1}^{\infty}$ a sofic approximation to $G$.  For \(\Phi \subset \Map(M, X)\), denote by \(N_\epsilon(\Phi; \rho_{X}^M | \rho_{Y}^M, \theta)\) the maximal cardinality of a \((\rho_{X}^M, \epsilon)\)-separated subset \(\Phi_0 \subset\Phi\) satisfying
		\[
		\rho_{Y}^M(\pi \circ \varphi, \pi \circ \psi) \leq \theta
		\]
		for all \(\varphi, \psi \in \Phi_0\). Define
		\begin{align*}
			h^\ep_{\Sigma}(\rho_X|\rho_Y; U, \delta | Y, \theta) &= \limsup_{i \to \infty} \frac{1}{\vol(M_i)} \log N_\epsilon \big( \Map(M_i, X, \rho_{X} \colon U, \delta); \rho_{X}^{M_i} | \rho_{Y}^{M_i} , \theta \big), \\
			h^\ep_{\Sigma}(\rho_X|\rho_Y, \theta) &= \inf_{U \subset G} \inf_{\delta > 0} h^\ep_{\Sigma}(\rho_X|\rho_Y; U, \delta | Y, \theta), \\
			h^\ep_{\Sigma}(\rho_X|\rho_Y) &= \inf_{\theta > 0} h^\ep_{\Sigma}(\rho_X|\rho_Y, \theta), \\
			h_{\Sigma}(\rho_X|\rho_Y)&= \sup_{\epsilon > 0} h^\ep_{\Sigma}(\rho_X|\rho_Y).
		\end{align*}
	\end{defn}
	\begin{lem}\label{X1}
		Let $\pi:(X,G)\to (Y,G)$ be a factor map between $G$-systems, and $\Sigma = (M_i)_{i=1}^{\infty}$  a sofic approximation to $G$.
		For a fixed continuous dynamically generating pseudometric \(\rho_Y\) on \(Y\), the quantity \(h_{\Sigma}(\rho_X|\rho_Y)\) does not depend on the choice of the metric \(\rho_X\) on \(X\).
	\end{lem}
	
	\begin{proof}
		Let \(\rho_1\) and \(\rho_2\) be compatible metrics on \(X\). Fix \(\varepsilon > 0\). Since the identity map \((X, \rho_1) \to (X, \rho_2)\) is uniformly continuous, there exists \(\varepsilon' > 0\) such that for any \(x_1, x_2 \in X\),
		\[
		\rho_2(x_1, x_2) \geq \varepsilon \quad \text{implies} \quad \rho_1(x_1, x_2) \geq \varepsilon'.
		\]
		
		Let $U_1\subset G$ be a precompact neighborhood of the identity and $\eta,\d_1>0$. Then by Lemma \ref{usefull}, for a sufficiently good sofic approximation $M$, there exists a precompact neighborhood of the identity $U_2\subset G$ and $\d_2>0$ such that $$\Map(M, X, \rho_2 \colon U_2,\d_2)\subset\Map(M, X, \rho_1 \colon U_1,\d_1).$$
		
		For any subset \(B \subset \Map(M, X, \rho_2 \colon U_2,\d_2)\), if \(B\) is \((\rho^{M}_{2}, \varepsilon)\)-separated, then it is also \((\rho^{M}_{1}, \varepsilon')\)-separated. Therefore, for any \(\theta > 0\),
		\[
		N_{\varepsilon}\left(\Map(M, X, \rho_2 \colon U_2,\d_2), \rho^{M}_{2}|\rho^M_{Y}, \theta\right) \leq N_{\varepsilon'}\left(\Map(M, X, \rho_1 \colon U_1,\d_1), \rho^{M}_{1}|\rho^M_{Y}, \theta\right).
		\]
		
		Taking appropriate limits, we obtain
		\[
		h_{\Sigma}(\rho_2|\rho_Y) \leq h_{\Sigma}(\rho_1|\rho_Y).
		\]
		
		By symmetry (exchanging the roles of \(\rho_1\) and \(\rho_2\)), the reverse inequality also holds. Hence,
		\[
		h_{\Sigma}(\rho_1|\rho_Y) = h_{\Sigma}(\rho_2|\rho_Y),
		\]
		which establishes the desired independence.	
	\end{proof}
	
	\begin{lem}\label{X2}
		Let $\pi:(X,G)\to (Y,G)$ be a factor map between $G$-systems, and $\Sigma = (M_i)_{i=1}^{\infty}$  a sofic approximation to $G$.
		For a fixed compatible metric \(\rho_X\) on \(X\), the quantity \(h_{\Sigma}(\rho_X|\rho_Y)\) does not depend on the choice of the metric \(\rho_Y\) on \(Y\).
	\end{lem}
	
	\begin{proof}
		Let \(\rho_1\) and \(\rho_2\) be compatible metrics on \(Y\). Fix \(\theta > 0\). Since the identity map \((Y, \rho_1) \to (Y, \rho_2)\) is uniformly continuous, there exists \(\theta' > 0\) such that for any \(y_1, y_2 \in Y\),
		\[
		\rho_1(y_1, y_2) \leq \theta' \quad \text{implies} \quad \rho_2(y_1, y_2) \leq \theta.
		\]
		
		For any \(\varepsilon > 0\), precompact neighborhood of the identity $U\subset G$, \(\delta > 0\), and sofic approximation $M$, we have
		\[
		N_{\varepsilon}\left(\Map(M, X, \rho_X \colon U,\d), \rho^{M}_{X}|\rho^M_{1}, \theta'\right) \leq N_{\varepsilon}\left(\Map(M, X, \rho_X \colon U,\d), \rho^{M}_{X}|\rho^M_{2}, \theta\right).
		\]
		
		Taking appropriate limits, we obtain
		\[
		h_{\Sigma}(\rho_{X}|\rho_{1}) \leq h_{\Sigma}(\rho_{X}|\rho_{2}).
		\]
		
		By symmetry (exchanging the roles of \(\rho_{1}\) and \(\rho_{2}\)), the reverse inequality also holds. Hence,
		\[
		h_{\Sigma}(\rho_{X}|\rho_{1})= h_{\Sigma}(\rho_{X}|\rho_{2}),
		\]
		which establishes the desired independence.
	\end{proof}
	By Lemma~\ref{X1} and \ref{X2}, we have the following result.
	
	\begin{thm}\label{main2}
		Let $\pi:(X,G)\to (Y,G)$ be a factor map between $G$-systems and $\Sigma = (M_i)_{i=1}^{\infty}$ a sofic approximation to $G$. Then $h_{\Sigma}(\rho_X|\rho_Y)$ does not depend on the choice of metrics $\rho_X$ and $\rho_Y$ on $X$ and $Y$ respectively.
	\end{thm}
	
	Next we will prove that the two definition methods of relative topological entropy are equivalent. Before that, we need the following Lemma (see \cite[Lemma 3.16]{luo17}).
	
	\begin{lem}\label{xy-y}
		Let \(\pi : X \to Y\) be a continuous map between two compact metrizable spaces \(X\) and \(Y\) with compatible metrics \(\rho_X\) and \(\rho_Y\) respectively. For any \(\delta > 0\) there exists an \(\epsilon > 0\) such that, for any \(x_1, x_2, \dots, x_m \in X\) with the condition that \(\rho_Y(\pi(x_i), \pi(x_j)) < \epsilon\) for all \(1 \leq i, j \leq m\), there exist \(z_1, z_2, \dots, z_m \in X\) such that \(\rho_X(x_i, z_i) < \delta\) for all \(1 \leq i \leq m\), and \(\pi(z_1) = \pi(z_2) = \cdots = \pi(z_m) \in Y\).
	\end{lem}
	
	\begin{prop}\label{eq}
		Let $\pi:(X,G)\to (Y,G)$ be a factor map between $G$-systems and $\Sigma = (M_i)_{i=1}^{\infty}$ a sofic approximation to $G$.	Let $\rho_X$ and $\rho_Y$ be two 1-bounded continuous dynamically generating pseudometrics \(\rho_X\) and \(\rho_Y\) on \(X\) and \(Y\) respectively. Then  we have
		\[
		h_{\Sigma}(\rho_X|Y) = h_{\Sigma}(\rho_X|\rho_Y).
		\]
		
		\begin{proof}
			
			The inequality $	h_{\Sigma}(\rho_X|Y) \leq h_{\Sigma}(\rho_X|\rho_Y)$ is immediate from the definitions. We now prove the reverse inequality.
			
			Fix \(\varepsilon > 0\), $U\subset G$ a precompact neighborhood of the identity and $\ep,\d_1>0$. Set \(\varepsilon_0 = 2\varepsilon\). It suffices to show that for sufficiently small parameters \(\delta' > 0\) and \(\theta > 0\), we have
			$$N_{\ep_0} \big( \Map(M, X, \rho_X \colon U,\d'); \rho^{M}_X | \rho^{M}_Y,\theta \big)\leq  N_{\ep}\big( \Map(M, X, \rho_X \colon U,\d); \rho^{M}_X | Y\big).$$

			Let \(\mathcal{F}\) be a maximal \((\rho^{M}_X, \varepsilon_0)\)-separated subset of $\Map(M_i, X, \rho_X \colon U,\d')$ such that \(\rho^{M}_Y(\pi \circ \varphi, \pi \circ \psi) \leq \theta\) for all \(\varphi, \psi \in \mathcal{F}\), and satisfying
			\[
			|\mathcal{F}| = N_{\ep_0} \big( \Map(M, X, \rho_X \colon U,\d'); \rho^{M}_X | \rho^{M}_Y,\theta \big).
			\]
			
			Choose \(\eta > 0\) such that \(\eta < \min\left\{\frac{\delta}{4}, \frac{\varepsilon}{2}\right\}\) and with the property that for any \(x_1, x_2 \in X\) with \(\rho_X(x_1, x_2) < \eta\), we have \(\rho_X(sx_1, sx_2) < \frac{\delta}{4}\) for all \(s \in U\). By Lemma~\ref{xy-y}, for sufficiently small \(\theta > 0\), there exists for each \(\varphi \in \mathcal{F}\) a map \(\widetilde{\varphi} : M \to X\) such that:
			\begin{itemize}
				\item For each fixed \(a \in M\), the value \(\pi(\widetilde{\varphi}(a))\) is independent of the choice of \(\varphi \in \mathcal{F}\);
				\item \(\rho_X(\widetilde{\varphi}(a), \varphi(a)) < \eta\) for all \(a \in M_i\) and all \(\varphi \in \mathcal{F}\).
			\end{itemize}
			Define \(\widetilde{\mathcal{F}} = \{\widetilde{\varphi} : \varphi \in \mathcal{F}\}\).
			
			\begin{claim}\label{c2}
				The set \(\widetilde{\mathcal{F}}\) is \((\rho^{M}_X, \varepsilon)\)-separated.
			\end{claim}
			\begin{proof}[Proof of Claim \ref{c2}]
				Let \(\widetilde{\varphi}_1, \widetilde{\varphi}_2 \in \widetilde{\mathcal{F}}\) be distinct. Since the original maps \(\varphi_1, \varphi_2 \in \mathcal{F}\) satisfy \( \rho^{M}_X(\varphi_1, \varphi_2) \geq \varepsilon_0 = 2\varepsilon\), we have
				\[
				\rho^{M}_X(\widetilde{\varphi}_1, \widetilde{\varphi}_2) \geq  \rho^{M}_X(\varphi_1, \varphi_2) -  \rho^{M}_X(\varphi_1, \widetilde{\varphi}_1) -  \rho^{M}_X(\widetilde{\varphi}_2, \varphi_2) > 2\varepsilon - \eta - \eta > \varepsilon.
				\]
				Hence, \(\widetilde{\mathcal{F}}\) is \(( \rho^{M}_X, \varepsilon)\)-separated.
			\end{proof}
			
			\begin{claim}\label{c3}
				We have \(\widetilde{\mathcal{F}} \subset \Map(M, X, \rho_X \colon U,\d)\).
			\end{claim}
			\begin{proof}[Proof of Claim \ref{c3}]
				Let \(\widetilde{\varphi} \in \widetilde{\mathcal{F}}\) with corresponding \(\varphi \in \mathcal{F}\). Since \(\rho^{M}_X(\varphi, \widetilde{\varphi}) < \eta\), the uniform continuity condition implies $\rho^{M}_X(s\varphi, s\widetilde{\varphi}) < \frac{\delta}{4}$ for all \(s \in U\).
				
				Now, for each \(s \in U\), we estimate:
				\[
				\rho^{M}_X(s\circ\widetilde{\varphi}, \widetilde{\varphi}\circ s) \leq \rho^{M}_X(s\circ\widetilde{\varphi}, s\circ\varphi) + \rho^{M}_X(s\circ\varphi, \varphi\circ s) + \rho^{M}_X( \varphi\circ s, \widetilde{\varphi}\circ s).
				\]
				The first term is bounded by \(\frac{\delta}{4}\) as above. The second term satisfies \(\rho^{M}_X(s\circ\varphi, \varphi\circ s) < \delta'\) since \(\varphi \in \Map(M, X, \rho_X \colon U,\d')\). For the third term,  we have
				\[
				\rho^{M}_X( \varphi\circ s, \widetilde{\varphi}\circ s) \leq \rho^{M}_X(\varphi, \widetilde{\varphi}) < \eta.
				\]
				Combining these estimates and using the conditions \(\delta' < \frac{\delta}{2}\) and \(\eta < \frac{\delta}{4}\), we obtain
				\[
				\rho^{M}_X(s\circ\widetilde{\varphi}, \widetilde{\varphi}\circ s) < \frac{\delta}{4} + \delta' + \eta < \delta.
				\]
				Therefore, \(\widetilde{\varphi} \in \Map(M, X, \rho_X \colon U,\d)\).
			\end{proof}
			
			From Claims \ref{c2} and \ref{c3}, it follows that
			\[
			N_{\varepsilon_0}\big(\Map(M, X, \rho_X \colon U,\d'),  \rho^{M}_X|\rho^{M}_Y, \theta\big) = |\mathcal{F}| = |\widetilde{\mathcal{F}}| \leq N_{\varepsilon}\big(\Map(M, X, \rho_X \colon U,\d),  \rho^{M}_X|Y\big).
			\]
			This establishes the desired inequality, completing the proof.
		\end{proof}
	\end{prop}
	
	It is interesting to define	relative sofic topological entropy of $(Y,G)$ with respect to the extension $(X,G)$.
	
	\begin{defn}
		Let $\pi:(X,G)\to (Y,G)$ be a factor map between $G$-systems. Let \(\rho_X\) and \(\rho_Y\) be continuous pseudometrics on \(X\) and \(Y\) respectively, $\Sigma = (M_i)_{i=1}^{\infty}$ a sofic approximation to $G$. Define
		\begin{align*}
			h^\ep_{\Sigma}(\rho_Y|\rho_X; U, \delta | X) &= \limsup_{i \to \infty} \frac{1}{\vol(M_i)} \log \Sep_\epsilon \big( \{\pi \circ \varphi : \varphi \in \Map(M_i, X, \rho_{X} \colon U, \delta)\}, \rho_{Y}^{M_i}\big), \\
			h^\ep_{\Sigma}(\rho_Y|\rho_X) &= \inf_{U \subset G} \inf_{\delta > 0} h^\ep_{\Sigma}(\rho_Y|\rho_X; U, \delta | X), \\
			h_{\Sigma}(\rho_Y|\rho_X)&= \sup_{\epsilon > 0} h^\ep_{\Sigma}(\rho_Y|\rho_X).
		\end{align*}
		
	\end{defn}
	\begin{lem}\label{Y1} 
		Let $\pi:(X,G)\to (Y,G)$ be a factor map between $G$-systems, and $\Sigma = (M_i)_{i=1}^{\infty}$  a sofic approximation to $G$.
		For a fixed continuous dynamically generating pseudometric \(\rho_Y\) on \(Y\), the quantity \(h_{\Sigma}(\rho_Y|\rho_X)\) does not depend on the choice of the continuous dynamically generating pseudometric \(\rho_X\) on \(X\).
	\end{lem}
	
	\begin{proof}
		Let \(\rho_1\) and \(\rho_2\) be two continuous dynamically generating pseudometrics on \(X\). Let $U_1\subset G$ be a precompact neighborhood of the identity and $\ep,\d_1>0$. Then  by Lemma \ref{usefull}, for a sufficiently good sofic approximation $M$, there exists a precompact neighborhood of the identity $U_2\subset G$ and $\d_2>0$ such that $$\Map(M, X, \rho_2 \colon U_2,\d_2)\subset\Map(M, X, \rho_1 \colon U_1,\d_1).$$
		Consequently, for any \(\varepsilon > 0\),
		\[
		h^{\varepsilon}_{\Sigma}(\rho_Y|\rho_2) \leq h^{\varepsilon}_{\Sigma}(\rho_Y|\rho_1),
		\]
		and thus
		\[
		h_{\Sigma}(\rho_Y|\rho_2) \leq h_{\Sigma}(\rho_Y|\rho_1).
		\]
		By symmetry (exchanging the roles of \(\rho_1\) and \(\rho_2\)), we also obtain the reverse inequality. Therefore,
		\[
		h_{\Sigma}(\rho_Y|\rho_1) = h_{\Sigma}(\rho_Y|\rho_2),
		\]
		which establishes the desired independence.
	\end{proof}
	
	\begin{lem}\label{Y2}
		Let $\pi:(X,G)\to (Y,G)$ be a factor map between $G$-systems, and $\Sigma = (M_i)_{i=1}^{\infty}$  a sofic approximation to $G$.
		For a fixed compatible metric \(\rho_X\) on \(X\), the quantity \(h_{\Sigma}(Y|X)\) does not depend on the choice of the continuous dynamically generating pseudometric \(\rho_Y\) on \(Y\).
	\end{lem}
	
	\begin{proof}
		Let \(\rho_1\) and \(\rho_2\) be two continuous dynamically generating pseudometrics on \(Y\). By symmetry, it suffices to show that for any \(\varepsilon > 0\), there exists \(\varepsilon' > 0\) such that for all  precompact neighborhood of the identity $U\subset G$, \(\delta > 0\),
		\begin{align*}
			&\Sep_{\varepsilon}\left(\{\pi \circ \varphi :\varphi \in \Map(M, X, \rho_X \colon U,\d)\}, \rho_1^M\right) \\&\leq \Sep_{\varepsilon'}\left(\{\pi \circ \varphi : \varphi \in \Map(M, X, \rho_X \colon U,\d)\},\rho_2^M\right).
		\end{align*}

		Since the identity map \((Y, \rho_1) \to (Y, \rho_2)\) is uniformly continuous, for any \(\varepsilon > 0\), there exists \(\varepsilon' > 0\) such that for all \(y_1, y_2 \in Y\),
		\[
		\rho_1(y_1, y_2) \geq \varepsilon \quad \text{implies} \quad \rho_2(y_1, y_2) \geq \varepsilon'.
		\]
		
		Now suppose \(\varphi, \psi \in \Map(M, X, \rho_X \colon U,\d)\) satisfy \(\rho_1^M(\pi \circ \varphi, \pi \circ \psi) \geq \varepsilon\). Then there exists \(p \in M\) such that
		\[
		\rho_1(\pi(\varphi(p)), \pi(\psi(p))) \geq \varepsilon,
		\]
		and consequently
		\[
		\rho_2(\pi(\varphi(p)), \pi(\psi(p))) \geq \varepsilon'.
		\]
		This implies that any \((\rho_2^M, \varepsilon')\)-separated subset of \(\{\pi \circ \varphi : \varphi \in \Map(M, X, \rho_X \colon U,\d)\}\) must have cardinality at least as large as that of any \((\rho_1^M, \varepsilon)\)-separated subset. Therefore,
		\begin{align*}
			&\Sep_{\varepsilon}\left(\{\pi \circ \varphi : \varphi \in \Map(M, X, \rho_X \colon U,\d)\}, \rho_1^M\right)\\& \leq \Sep_{\varepsilon'}\left(\{\pi \circ \varphi :\varphi \in \Map(M, X, \rho_X \colon U,\d)\},\rho_2^M\right),
		\end{align*}
		as required.
	\end{proof}

	By Lemma~\ref{Y1} and \ref{Y2}, we have the following result.

	\begin{thm}\label{mainY}
		Let $\pi:(X,G)\to (Y,G)$ be a factor map between $G$-systems and $\Sigma = (M_i)_{i=1}^{\infty}$ a sofic approximation to $G$. Then $h_{\Sigma}(\rho_Y|\rho_X)$ does not depend on the choice of metrics $\rho_X$ and $\rho_Y$ on $X$ and $Y$ respectively.
	\end{thm}
	
	We shall denote this common value by $h_\Sigma(Y|X)$, and call it the \emph{relative sofic topological entropy} of $(Y,G)$  with respect to the extension $(X,G)$.

	\begin{prop}
		Let $\pi:(X,G)\to (Y,G)$ be a factor map between $G$-systems, and $\Sigma = (M_i)_{i=1}^{\infty}$  a sofic approximation to $G$.  Then
		\begin{enumerate}
			\item \(h_{\Sigma}(Y|X) \leq h_{\Sigma}(X)\),
			\item \(h_{\Sigma}(Y|X) \leq h_{\Sigma}(Y)\),
			\item \(h_{\Sigma}(X|Y) \leq h_{\Sigma}(X)\).
		\end{enumerate}
		
		\begin{proof}
			Let \(\rho_X\) and \(\rho_Y\) be compatible metrics on \(X\) and \(Y\) respectively. Define a new metric \(\rho'_X\) on \(X\) by
			\[
			\rho'_X(x_1, x_2) = \rho_X(x_1, x_2) + \rho_Y(\pi(x_1), \pi(x_2))
			\]
			for any \(x_1, x_2 \in X\). Since \(\rho'_X\) is compatible with the topology on \(X\), we may replace \(\rho_X\) with \(\rho'_X\) and assume without loss of generality that
			\[
			\rho_X(x_1, x_2) \geq \rho_Y(\pi(x_1), \pi(x_2))
			\]
			for all \(x_1, x_2 \in X\).
			
			First, we prove \(h_{\Sigma}(Y|X) \leq h_{\Sigma}(X)\). 
			This follows from the observation that for any \(\varepsilon > 0\), \(\delta > 0\), $U\subset G$ a precompact neighborhood of the identity, we have
			\[
			\Sep_\varepsilon\big(\{\pi \circ \varphi : \varphi \in \Map(M, X, \rho_X \colon U,\d)\}, \rho^{M}_{Y}\big) \leq \Sep_\varepsilon\big(\Map(M, X, \rho_X \colon U,\d), \rho^{M}_{X}\big).
			\]
			
			Next, we prove \(h_{\Sigma}(Y|X) \leq h_{\Sigma}(Y)\) and it suffices to show the follow claim.
			\begin{claim}
				For any \(\delta > 0\), $U\subset G$ a precompact neighborhood of the identity, there exists $\d'>0$, such that
				\[
				\{\pi \circ \varphi : \varphi \in \Map(M, X, \rho_X \colon U,\d)\} \subset \Map(M, X, \rho_X \colon U,\d).
				\]
			\end{claim}

			\begin{proof}
				Since $\pi$ is continuous and $X$ is compact, we can find $\delta' > 0$ small enough such
				that $\delta'(\diam(Y, \rho_Y))^2 \leq \delta^2/4$ and for any $x, x' \in X$ with $\rho_X(x, x') \leq \sqrt{\delta'}$ one has
				$\rho_Y(\pi(x), \pi(x')) \leq \delta/2$.
				
				Fix	$s\in U$. Let	$A := \left\{ p \in M : \rho_X( s \circ \varphi(p), \varphi(s.p))\leq \sqrt{\d'} \right\}$. By Markov’s inequality one
				has $\ov{\vol}(M\setminus A)\leq \sqrt{\d'}$. For each $a \in A$, by the choice of $\delta'$ one has
				\[
				\rho_Y(s\circ\pi(\varphi(a)), \pi(\varphi(s.a))) = \rho_Y(\pi(s\varphi(a)), \pi(\varphi(s.a))) \leq \delta/2.
				\]
				
				Thus
				\begin{align*}
					\rho_Y^M(s\circ \pi \circ \varphi, \pi \circ \varphi \circ s)
					&\leq \ov{\vol}(A) \cdot \frac{\delta}{2} + \ov{\vol}(M\setminus A) \cdot \diam(Y, \rho_Y) \\
					&\leq \frac{\delta}{2} + \sqrt{\d'}(\diam(Y, \rho_Y))^2 \\
					&\leq \frac{\delta}{2} + \frac{\delta}{2}=\delta 
					,
				\end{align*}
				and hence $\pi \circ \varphi \in \Map(M, X, \rho_X \colon U,\d)$.
			\end{proof}
			
			Finally, the inequality \(h_{\Sigma}(X|Y) \leq h_{\Sigma}(X)\) is immediate from Definition.
		\end{proof}
	\end{prop}
	
	Now we can prove the additive inequality.
	
	\begin{thm}\label{add}
		Let $\pi:(X,G)\to (Y,G)$ be a factor map between $G$-systems, and $\Sigma = (M_i)_{i=1}^{\infty}$  a sofic approximation to $G$.  Then
		\[
		h_{\Sigma}(X) \leq h_{\Sigma}(X|Y) + h_{\Sigma}(Y|X).
		\]
		
		\begin{proof}
			Let \(\rho_X\) and \(\rho_Y\) be compatible metrics on \(X\) and \(Y\) respectively. By Proposition~\ref{eq}, Theorem~\ref{main2} and \ref{mainY}, it suffices to show that
			\[
			h_{\Sigma}(X) \leq h_{\Sigma}(\rho_X|\rho_Y) + h_{\Sigma}(\rho_Y|\rho_X).
			\]
			
			Fix \(\varepsilon, \varepsilon', \delta > 0\) and $U\subset G$ a precompact neighborhood of the identity.  Let \(\varphi\) be a maximal \((\rho_{X}^{M}, \varepsilon)\)-separated subset of $\Map(M, X, \rho_X \colon U,\d)$. Let \(\Psi\) be a maximal \((\rho_{Y}^{M}, \varepsilon')\)-separated subset of \(\{\pi \circ \varphi :\varphi \in \Map(M, X, \rho_X \colon U,\d)\}\).
			
			For each \(\varphi \in \Phi\), there exists \(\psi_{\varphi} \in \Psi\) such that \(\rho_{Y}^{M_i}(\pi \circ \varphi, \psi_{\varphi}) < \varepsilon'\). By the pigeonhole principle, there exists a subset \(\Phi_0 \subset \Phi\) with
			\[
			|\Phi_0| \geq \frac{|\Phi|}{|\Psi|}
			\]
			and a fixed \(\psi_0 \in \Psi\) such that \(\psi_{\varphi} = \psi_0\) for all \(\varphi \in \Phi_0\).
			
			For any \(\varphi, \varphi' \in \Phi_0\), we have
			\[
			\rho_{Y}^{M}(\pi \circ \varphi, \pi \circ \varphi') \leq \rho_{Y}^{M}(\pi \circ \varphi, \psi_0) + \rho_{Y}^{M}(\psi_0, \pi \circ \varphi') < 2\varepsilon'.
			\]
			Since \(\Phi_0\) is contained in $\Phi$ then we have 
			\[
			|\Phi_0| \leq N_{\varepsilon}\big(\Map(M, X, \rho_X \colon U,\d), \rho_{X}^{M}|\rho_{Y}^{M}, 2\varepsilon'\big).
			\]
			
			Combining the inequalities, we obtain
			\begin{align*}
				|\Phi| &\leq |\Psi| \cdot |\Phi_0| \\&\leq N_{\varepsilon'}\big(\{\pi \circ \varphi :\varphi \in \Map(M, X, \rho_X \colon U,\d)\}, \rho_{Y}^{M}\big) \cdot N_{\varepsilon}\big(\Map(M, X, \rho_X \colon U,\d), \rho_{X}^{M}|\rho_{Y}^{M}, 2\varepsilon'\big).
			\end{align*}

			Taking logarithms and appropriate limits, we get
			\[
			h^{\varepsilon}_{\Sigma}(X) \leq h^{\varepsilon}_{\Sigma}(\rho_X; U, \delta) \leq h^{\varepsilon}_{\Sigma}(\rho_X|\rho_Y; U, \delta | Y, 2\ep') + h^{\varepsilon'}_{\Sigma}(\rho_Y; U, \delta|X).
			\]
			
			Since $U\subset G$  and \(\delta > 0\) are arbitrary, we have
			\[
			h^{\varepsilon}_{\Sigma}(X) \leq h^{\varepsilon}_{\Sigma}(\rho_X|\rho_Y, 2\varepsilon') + h_{\Sigma}(\rho_Y|\rho_X).
			\]
			
			Taking the limit as \(\varepsilon' \to 0\), we obtain
			\[
			h^{\varepsilon}_{\Sigma}(X) \leq h^{\varepsilon}_{\Sigma}(\rho_X|\rho_Y) + h_{\Sigma}(\rho_Y|\rho_X) \leq h_{\Sigma}(\rho_X|\rho_Y) + h_{\Sigma}(\rho_Y|\rho_X).
			\]
			
			Finally, taking the limit as \(\varepsilon \to 0\), we conclude
			\[
			h_{\Sigma}(X) \leq h_{\Sigma}(X|Y) + h_{\Sigma}(Y|X),
			\]
			which completes the proof.
		\end{proof}
	\end{thm}

	
	%
	
	\section{Relative sofic measure entropy and variational principle}\label{sec4}
	In this section, we study relative sofic measure entropy and show the relative version of the variational principle.

	Denote by $\pr(X)$ the set of all Borel probability measures on $X$. We say that $\mu \in \pr(X)$ is \textit{$G$-invariant} if $\mu = g\mu:=\mu\circ g^{-1}$ for each $g \in G$. Denote by $\mathrm{Prob}(X,G)$  the set of all $G$-invariant  measures on $X$. Note that both $\pr(X)$ and $\pr(X,G)$ are convex compact metric spaces when they are endowed with the weak$^*$-topology. The support of of a measure $\mu \in \pr(X)$, denoted by $\supp\mu$, is defined as follow:
	\begin{align*}
		\supp\mu= \{x\in X: \text{for any open neighbourhood}\ U\ \text{of}\ x, \ \mu(U)>0\}.
	\end{align*}
	Equivalently, $\supp\mu$ is the smallest closed subset $C$ of $X$ that satisfies $\mu(C)=1$.
	\begin{defn}
		Let \( U \subset G \) be a precompact neighborhood of identity, \( O \subset \pr(X) \) an open neighborhood of \( \mu \), \( \epsilon > 0 \). Suppose that \( M \) is a local \( G \)-space with finite measure \( \mathrm{vol} \). As before we let \( \mathrm{vol} \) denote the normalized volume on \( M \), so that \( \mathrm{vol}(M') = \frac{\mathrm{vol}(M')}{\mathrm{vol}(M)} \) for any \( M' \subset M \). Denote by \( \mathrm{Map}(M, X, \rho: U, \delta, O) \) the set of all maps \( \varphi \in \mathrm{Map}(M, X, \rho: U, \delta) \) such that \( \varphi_* \ov{\mathrm{vol}} \in O \) (where $\varphi_* \ov{\mathrm{vol}}$ is the measure on defined by $\varphi_* \ov{\mathrm{vol}}(E)=\ov{\mathrm{vol}}(\varphi^{-1}(E))$).
	\end{defn}
	
	\begin{defn}
		Let $\pi:(X,G)\to (Y,G)$ be a factor map between $G$-systems. Let \(\rho_X\) and \(\rho_Y\) be continuous pseudometrics on \(X\) and \(Y\) respectively, $\Sigma = (M_i)_{i=1}^{\infty}$ a sofic approximation to $G$. Define
		\begin{align*}
			h_{\Sigma,\mu}^\ep(\rho_X|Y) &=\inf_{O} \inf_{U \subset G} \inf_{\delta > 0} \limsup_{i \to \infty} \frac{1}{\vol(M_i)} \log N_\epsilon \big( \Map(M_i, X, \rho_{X} \colon U,\d,O); \rho_{X}^{M_i} | Y\big), \\
			h_{\Sigma,\mu}(\rho_X|Y) &= \sup_{\epsilon > 0} h_{\Sigma,\mu}^\ep(\rho_X|Y).
		\end{align*}
		where the first infimum  is over all open neighborhoods $O$ of $\mu$ in $\pr(X)$ and the second infimum is over all precompact neighborhoods $U$ of the identity in $G$.
		When $(Y,G)$ is trivial, we recover standard sofic measure entropy:
		$$	h_{\Sigma,\mu}(\rho_X):=h_{\Sigma,\mu}(\rho_X|Y).$$
	\end{defn}
	
	\begin{thm}
		Let $\Sigma = (M_i)_{i=1}^{\infty}$ be a sofic approximation to $G$,  $\rho_1$ and $\rho_2$ be two bounded dynamically generating continuous pseudometrics on $X$. Let $\pi:(X,G)\to (Y,G)$ be a factor map between $G$-systems, and $\mu\in\pr(X,G)$. Then we have 
		$$h_{\Sigma,\mu}(\rho_1|Y)=h_{\Sigma,\mu}(\rho_2|Y).$$
	\end{thm}
	\begin{proof}
		The proof is nearly identical to the proof of Theorem \ref{main1}. We omit it.
	\end{proof}
	We shall denote this common value by $h_{\Sigma,\mu}(X|Y)$, and call it the \emph{relative sofic measure entropy} of $(X,G)$ with respect to the factor $(Y,G)$.
	
	For future applications, it is nice to have a slightly more general version. If
	$A \subset \pr(X)$ is any set, then define the entropy
	\[
	h_{\Sigma, A}(X|Y) = \lim_{\epsilon \searrow 0} 
	\inf_{\substack{O \supset A \\ \text{open}}} 
	\inf_{\substack{U \subset G \\ \text{precompact}}} 
	\inf_{\delta > 0} 
	\limsup_{n \to \infty} 
	\mathrm{vol}(M_n)^{-1} \log N_{\epsilon}\bigl(\Map(M_n, X, \rho : U, \delta, O),\rho^{M_n}|Y\bigr),
	\]
	where $O$ varies over all open sets containing $A$ and $U$ varies over all precompact
	open subsets $U$ of $G$.
	
	This definition is the same as the definition of measure entropy, with the exception that the open set $O$ is required to contain $A$ instead of $\mu$. In particular, if
	$A = \{\mu\}$ then $h_{\Sigma, A}(X|Y) = h_{\Sigma, \mu}(X|Y)$; if
	$A = \pr(X)$ then $h_{\Sigma, A}(X|Y) = h_{\Sigma}(X|Y)$. 
	
	Next we will prove the variational principle and	we need some preparation. We need the following Lemma (see \cite[Lemma 7.3]{B25}).
	\begin{lem}\label{G-inv}
		Let $(X,G)$ be a $G$-system and $O$ an open neighborhood of $\pr(X,G)$ in $\pr(X)$. Then there exist a precompact open set $U \subset G$ and $\delta > 0$ such that if $M$ is a $(U, \delta)$-sofic approximation to $G$ and $\varphi \in \mathrm{Map}(M, X, \rho: U, \delta)$ then $\varphi_* \overline{\mathrm{vol}} \in O$.
	\end{lem}
	
	The following Theorem is a  classical result, for a proof see \cite{W75}.
	
	\begin{thm}\label{port}
		Let $X$ be a compact metric space, and let $\{\mu_n\}_{n=1}^\infty$ and $\mu$ be Borel probability measures on $X$. The following conditions are equivalent:
		
		\begin{enumerate}
			\item $\mu_n \to \mu$ in the weak$^*$-topology;
			
			\item For every open set $U \subseteq X$,	$\liminf_{n\to\infty} \mu_n(U) \geq \mu(U)$;
			
			\item For every closed set $F \subseteq X$,
			$
			\limsup_{n\to\infty} \mu_n(F) \leq \mu(F)
			$;
			
			\item For every Borel set $A \subseteq X$ with $\mu(\partial A) = 0$ (where $\partial A$ denotes the topological boundary of $A$), $	\lim_{n\to\infty} \mu_n(A) = \mu(A).$

		\end{enumerate}
	\end{thm}
	
	\begin{thm}\label{vxy}
		Let $\pi:(X,G)\to (Y,G)$ be a factor map between $G$-systems. Then for any sofic approximation $\Sigma$ to $G$ and closed subset $A \subset \pr(X)$,
		\[
		h_{\Sigma, A}(X|Y) = \sup_{\mu \in A \cap \pr(X,G)}h_{\Sigma, \mu}(X|Y).
		\]
	\end{thm}
	
	\begin{proof}
		Let \(\rho_X\) be a compatible metrics on \(X\). Firstly, we show that	\[
		\sup_{\mu \in A \cap \pr(X,G)}h_{\Sigma, \mu}(\rho_X|Y)\ge h_{\Sigma, A}(\rho_X|Y).
		\]
		Without loss of generality, one can assume that
		$h_{\Sigma, A}(\rho_X|Y) \neq -\infty$. Fix $\epsilon > 0$. For any $\d>0$, $k\in\N$, $U\subset G$ a precompact neighborhood of the identity and an open
		neighborhood $O$ of $A$ in $\pr(X)$, let $d_k(U, \delta, O)$	be a maximal $(\rho_X^{M_k}, \epsilon)$-separated subset of $\Map(M_k, X, \rho_X : U, \delta, O)$ satisfying $\pi \circ \varphi=\pi \circ \psi$
		for all \(\varphi, \psi \in 	d_k(U, \delta, O)\). Define
		\[
		h^{\epsilon}_{\Sigma, A}(\rho_X|Y) = \inf_{O} \inf_{U} \inf_{\delta > 0} 
		\limsup_{k \to \infty} \frac{1}{\mathrm{vol}(M_k)}\log |d_k(U, \delta, O)|.
		\]
		By the definition of relative sofic entropy, we have
		\[
		h_{\Sigma, A}(\rho_X|Y) = \lim_{\epsilon \searrow 0} h^{\epsilon}_{\Sigma, A}(\rho_X|Y).
		\]
		By the assumption of $h_{\Sigma, A}(\rho_X|Y) \neq -\infty$, one can see that for all $U, \delta, O$ the set $d_k(U, \delta, O)$
		is nonempty for infinitely many $k$.
		
		For any $\varphi \in \Map(M_k, X)$, define 
		\[
		\delta_{\varphi_*\ov{\mathrm{vol}}}(B) = 
		\begin{cases} 
			1 & \text{if } \varphi_*\ov{\mathrm{vol}} \in B, \\
			0 & \text{otherwise}.
		\end{cases}
		\]
		Then $\delta_{\varphi_*\ov{\mathrm{vol}}} \in \pr(\pr(X))$.
		Take $\gamma_k(U, \delta, O) = |d_k(U, \delta, O)|^{-1} \sum_{\varphi \in d_k(U, \delta, O)} \delta_{\varphi_*\ov{\mathrm{vol}}}$ whenever $d_k(U, \delta, O)$ is nonempty. Then $\gamma_k(U, \delta, O) \in \pr(\pr(X))$.
		
		For any $U, \delta, O$, one can take an increasing sequence $\{n_k\}_{k=1}^{\infty}$ such that
		\begin{equation}\label{star}
			\limsup_{k \to \infty} \frac{1}{\mathrm{vol}(M_k)}\log |d_k(U, \delta, O)| = \lim_{k \to \infty} \frac{1}{\mathrm{vol}(M_{n_k})} \log |d_{n_k}(U, \delta, O)|. 
		\end{equation}

		It is well known that $\pr(X)$  is a compact metric space in the
		weak* topology when $X$ is a compact metric space(see \cite{W75}). Then $\pr(\pr(X))$  is also a compact metric space in the
		weak* topology. Therefore, without loss of generality, we can assume	$ \gamma_{n_k}(U, \delta, O))$ conveges to some measure $\gamma(U, \delta, O) \in \pr(\pr(X))$ when $k \to \infty$. Let $\Gamma$ be an accumulation point
		of $\gamma(U, \delta, O)$ as $U \nearrow G$, $\delta \searrow 0$ and $O \searrow A$.
		
		By Lemma \ref{G-inv} and the outer regularity of probability measure, $\Gamma(\pr(X,G)) = 1$. By the definition of $\Map(M_k, X, \rho_X : U, \delta, O)$, $\Gamma(A) = 1$. As $\pr(X,G)$ and $A$ are closed subset of $\pr(\pr(X))$, $\pr(X,G)\cap A$ contains $\supp \Gamma$. Take $\mu\in \supp \Gamma$. Now we show that $h^{\epsilon}_{\Sigma, \mu}(\rho_X|Y) \geq h^{\epsilon}_{\Sigma, A}(\rho_X|Y)$.
		
		Let $W_{\mu}$ be an open neighborhood of $\mu$ in $\pr(X)$. Then $\Gamma(W_{\mu}) > 0$. Fix
		$U_0, \delta_0, O_0$. By Theorem~\ref{port}, there exist
		$U \supset U_0$, $0 < \delta \leq \delta_0$ and  an  open neighborhood $W \subset O_0$ of $A$ such that	$\gamma(U, \delta, W)(W_{\mu}) \geq \Gamma(W_{\mu})/2 > 0$.
		For  $U, \delta, W$, take the increasing sequence $\{n_k\}_{k=1}^{\infty}$ mentioned above. By  Theorem~\ref{port} again, there exists $N$ such that if $k > N$ then $\gamma_{n_k}(U, \delta, W)(W_{\mu}) \geq \Gamma(W_{\mu})/4 > 0$. Note that
		\[
		\{\varphi \in d_{n_k}(U, \delta, W) : \varphi_*\ov{\mathrm{vol}}\in W_{\mu}\} = d_{n_k}(U, \delta, W) \cap \Map(M_{n_k}, X, \rho_X : U, \delta, W \cap W_{\mu}).
		\]
		Therefore,
		\[
		\frac{|d_{n_k}(U, \delta, W) \cap \Map(M_{n_k}, X, \rho_X : U, \delta, W \cap W_{\mu})|}{|d_{n_k}(U,\delta, W)|} \geq \frac{\Gamma(W_{\mu})}{4}.
		\]
		Thus
		\begin{align*}
			N_{\epsilon}\bigl(\Map(M_{n_k}, X, \rho_X : U, \delta, W_{\mu}),\rho^{M_{n_k}}|Y\bigr) 
			&\geq N_{\epsilon}\bigl(\Map(M_{n_k}, X, \rho_X : U, \delta, W\cap W_{\mu}),\rho_X^{M_{n_k}}|Y\bigr) \\
			&\geq |d_{n_k}(U, \delta, W)|\Gamma(W_{\mu})/4.
		\end{align*}
		Take the logarithm of both sides and divide by  $\mathrm{vol}(M_{n_k})$. Then take the limit superior as $k \to \infty$, and finally take the infimum over the parameters $U$, $\delta$, and $W$. Then using the (\ref{star}), we have
		\[
		\inf_{U \subset G} \inf_{\delta > 0} \limsup_{k \to \infty} \frac{1}{\mathrm{vol}(M_k)} \log N_{\epsilon}\bigl(\Map(M_k, X, \rho_X : U, \delta, W_{\mu}),\rho_X^{M_k}|Y\bigr) \geq h^{\epsilon}_{\Sigma, A}(\rho_X|Y).
		\]
		Take the infimum over  all open neighborhoods $W_{\mu}$ of $\mu$ in $\pr(X)$ . Then we have $h^{\epsilon}_{\Sigma, \mu}(\rho_X|Y) \geq h^{\epsilon}_{\Sigma, A}(\rho_X|Y)$.
		Thus
		\[
		\sup_{\mu \in A \cap \pr(X,G)}h_{\Sigma, \mu}(\rho_X|Y)\ge 	h_{\Sigma, A}(\rho_X|Y).
		\]
		The opposite inequality
		\[
		\sup_{\mu \in A \cap \pr(X,G)}h_{\Sigma, \mu}(\rho_X|Y)	\leq h_{\Sigma, A}(\rho_X|Y) 
		\]
		is immediate. Therefore the proof of variational principle is complete.
	\end{proof}

	\section{Topological pressure and the variational principle}\label{sec5}
	In this section, we will deﬁne the sofic topological pressure and establish some basic properties	of it.

	\subsection{Variational principle}
	
	Let $\Sigma = (M_i)_{i=1}^{\infty}$ be a sofic approximation to $G$ and $(X,G)$ be a $G$-system. Let $\Omega$ be a real-valued continuous function on $X$, $\rho$ a continuous pseudometric on $X$.
	\begin{defn}
		Let $U \subset G$ be a precompact neighborhood of the identity, and $\delta,\ep > 0$.  Suppose that \( M \) is a local \( G \)-space with finite measure \( \mathrm{vol} \).	 We define
		\[
		M_{\Sigma}^{\ep}(\Omega, M, X, \rho: U, \delta) = \sup_{\mathcal{R}} \sum_{\varphi \in \mathcal{R}} \exp\left( \int_M \Omega(\varphi(p)) \, d{\mathrm{vol}}(p) \right),
		\]
		where $\mathcal{R}$ runs over $(\rho^{M}, \varepsilon)$-separated subsets of $\mathrm{Map}(M, X, \rho: U, \delta)$. Of course, the value of the right hand side does not change if $\mathcal{R}$  runs over maximal $(\rho^{M}, \varepsilon)$-separated subsets of $\mathrm{Map}(M, X, \rho: U, \delta)$.
	\end{defn}
	
	Now we define the sofic topological pressure of $\Omega$.
	
	\begin{defn}
		Let $\Sigma = (M_i)_{i=1}^{\infty}$ be a sofic approximation to $G$ and $(X,G)$ be a $G$-system.	We define
		\begin{align*}
			P_{\Sigma}^{\ep}(\Omega, X, G, \rho, U, \delta) &= \limsup_{i \to \infty} \frac{1}{\vol(M_i)}  \log M_{\Sigma}^{\ep}(\Omega, M_i, X, \rho: U, \delta), \\
			P_{\Sigma}^{\ep}(\Omega, X, G, \rho, U) &= \inf_{\delta > 0} P_{\Sigma}^{\ep}(\Omega, X, G, \rho, U, \delta), \\
			P_{\Sigma}^{\ep}(\Omega, X, G, \rho) &= \inf_{U} P_{\Sigma}^{\ep}(\Omega, X, G, \rho, U), \\
			P_{\Sigma}(\Omega, X, G, \rho) &= \sup_{\varepsilon > 0} P_{\Sigma}^{\ep}(\Omega, X, G, \rho),
		\end{align*}
		where $U$ in the third line is over all precompact neighborhoods $U$ of the identity in $G$. 
		
		If $\mathrm{Map}(M, X, \rho: U, \delta)=\emptyset$ for all large enough $i$, we set $	P_{\Sigma}^{\ep}(\Omega, X, G, \rho, U, \delta) = -\infty$.
	\end{defn}
	\begin{rem}\label{rem1}
		Actually, the deﬁnitions are straightforward
		generalizations of the countable case (see for example~\cite{C13}).	When $\Omega=0$, $P_{\Sigma}(0, X, G, \rho)$ is the standard sofic topological entropy $h_{\Sigma}(\rho)$ in Definition~\ref{def1}.
	\end{rem}

	The following result illustrates that the deﬁnition of soﬁc topological pressure does not depend on the
	choice of dynamically generating continuous pseudometrics
	
	\begin{thm}\label{main4}
		Let $\Sigma = (M_i)_{i=1}^{\infty}$ be a sofic approximation to $G$, $(X,G)$ be a $G$-system and $\Omega$ be a real-valued continuous function on $X$. Let $\rho_1$ and $\rho_2$ be two 1-bounded dynamically generating continuous pseudometrics on $X$. Then we have 
		$$	P_{\Sigma}(\Omega, X, G, \rho_1)=	P_{\Sigma}(\Omega, X, G, \rho_2)$$
	\end{thm}
	\begin{proof}
		The proof is nearly identical to the proof of Theorem~\ref{main1}. We omit it.
	\end{proof}
	
	Since $\pr(X)$	is endowed with the weak$^*$-topology, one can use the following method to define sofic measure entropy. We denote by $C(X)$ the set of all real-valued continuous function on $X$.
	\begin{defn}
		Let $\Sigma = (M_i)_{i=1}^{\infty}$ be a sofic approximation to $G$, $(X,G)$ be a $G$-system and $\mu\in\pr(X,G)$.	Let $U \subset G$ be a precompact neighborhood of the identity, $F \subset C(X)$ be a finite set, and $\delta > 0$.  Suppose that \( M \) is a local \( G \)-space with finite measure \( \mathrm{vol} \). Denote by \( \mathrm{Map}_\mu(M, X, \rho: U, F,\delta) \) the set of all maps \( \varphi \in \mathrm{Map}(M, X, \rho: U, \delta) \) such that
		$$ \left| \int f(\varphi(p)) \, d\ov{\mathrm{vol}}(p) - \int f \, d\mu \right| < \delta$$ for all $f \in F$.
		Define
		\begin{align*}
			h_{\Sigma,\mu}^\ep(\rho:U,F,\d) &=   \limsup_{i \to \infty} \frac{1}{\vol(M_i)} \log \sep_\epsilon \big( \Map_\mu(M_i, X, \rho \colon U,F,\d)\big), \\
			h_{\Sigma,\mu}^\ep(\rho:U,F)&=\inf_{\delta > 0}h_{\Sigma,\mu}^\ep(\rho:U,F,\d),\\
			h_{\Sigma,\mu}^\ep(\rho:U)&=\inf_{F\subset C(X)}h_{\Sigma,\mu}^\ep(\rho:U,F),\\
			h_{\Sigma,\mu}^\ep(\rho)&=\inf_{U \subset G}h_{\Sigma,\mu}^\ep(\rho:U),\\
			h_{\Sigma,\mu}(\rho) &= \sup_{\epsilon > 0} h_{\Sigma,\mu}^\ep(\rho).
		\end{align*}
	\end{defn}

	For locally compact soﬁc group actions, we also have the variational principle about soﬁc topological pressure.
	
	\begin{thm}\label{main3}
		Let $\Sigma = (M_i)_{i=1}^{\infty}$ be a sofic\ approximation sequence for $G$, $(X,G)$ be a $G$-system and $h$ be a real-valued continuous function on $X$. Then
		\[
		\
		P_\Sigma(\Omega, X, G) = \sup \left\{ h_{\Sigma,\mu}(X) + \int_X \Omega\  d\mu : \mu \in \pr
		(X,G) \right\}.
		\]
		In particular, if $P_\Sigma(\Omega, X, G) \neq -\infty$ then $\pr(X,G)$ 
		is non-empty.
	\end{thm}
	
	We need the following lemma which was proved by  Singh in (\cite[Theorem 4.2.1]{S16}) for the case $\Omega = 0$.
	We adapt the argument there to deal with general functions $\Omega\in C(X)$.
	
	\begin{lem}\label{lma3}
		Let $\Sigma = (M_i)_{i=1}^{\infty}$ be a sofic\ approximation sequence for $G$, $(X,G)$ be a $G$-system and $\Omega$ be a real-valued continuous function on $X$. Then
		\[
		\
		P_\Sigma(\Omega, X, G)\leq\sup \left\{ h_{\Sigma,\mu}(X) + \int_X \Omega \ d\mu : \mu \in \pr
		(X,G) \right\}.
		\]
	\end{lem}

	\begin{proof}
		%
		Let $e \in U_1 \subset U_2 \subset \cdots \subset U_n \subset \cdots$ be an increasing sequence of precompact neighborhoods of the identity in $G$. Since $X$ is a compact metric space, there exists a sequence $\{f_m\}_{m \in \mathbb{N}}$ in $C(X)$ such that
		$\{f_m\}_{m \in \mathbb{N}}$ is dense in $C(X)$. Let $n \in \mathbb{N}$ and $F_n = \{1,\Omega, f_1, \dots, f_n\}$. Using the continuity of the action and compactness of $X$, choose finite subsets $e \in A_n \subset U_n$ such that for every $g \in U_n$ and $f \in F_n$, there exists some $a_g \in A_n$ with
		\[
		\| f\circ g - f\circ a_g  \|_\infty \leq \frac{1}{4n}.
		\]
		Moreover, choose $0 < \delta_n < \left( \frac{1}{n \sup_{f \in F_n} \| f \|_\infty} \right)^2$ small enough so that for every $f \in F_n$ and any $x, y \in X$, $\rho(x, y) \leq \delta_n$ implies $|f(x) - f(y)| \leq 1/n$.
		
		Since the space $\pr(X)$ of all Borel probability measures on $X$ is weak* compact, there exists a finite set $B \subset \pr(X)$ such that for every $\varphi \in \mathrm{Map}(M, X, \rho: U_n, \delta_n)$, there is a measure $\mu_\varphi \in B$ satisfying
		\begin{equation} \label{eq:4.1}
			\left| \int f\circ \gamma\  d\mu_\varphi - \int f\circ \gamma \ d(\varphi_*\ov{\mathrm{vol}}) \right| < \frac{1}{4n}
		\end{equation}
		for every $f \in F_n$ and $\gamma \in A_n$.
		
		Let $S_n \subset \mathrm{Map}(M, X, \rho: U_n, \delta_n)$ be a maximal $(\rho^{M}, \varepsilon)$-separated set  such that
		\[
		M_{\Sigma}^{\ep}(\Omega, M, X, \rho: U_n, \delta_n) \leq \exp(1) \cdot \sum_{\varphi \in S_n} \exp\left( \int_M \Omega(\varphi(p)) \, d{\mathrm{vol}}(p) \right).
		\]
		By the pigeonhole principle, there exists a measure $\mu_n \in B$ such that
		
		\[
		|B| \cdot \sum_{\varphi \in \{ \varphi \in S_n : \mu_\varphi = \mu_n \}} \exp\left( \int_M \Omega(\varphi(p)) \, d{\mathrm{vol}}(p) \right) \geq \sum_{\varphi \in S_n} \exp\left( \int_M \Omega(\varphi(p)) \, d{\mathrm{vol}}(p) \right).
		\]
		
		%
		
		Denote $R(M,n)=\{ \varphi \in S_n : \mu_\varphi = \mu_n \}$ and let $\gamma = e \in G$ in \eqref{eq:4.1}, we obtain
		
		\[
		\left| \int \Omega\  d\mu_\varphi - \int_M \Omega(\varphi(p)) \, d{\ov{\mathrm{vol}}}(p) \right| < \frac{1}{4n}\quad \text{for
			all } \varphi \in R(M, n).
		\] Then we have
		
		\[
		\exp\left( \vol (M)\cdot\int \Omega \, d\mu_n + \frac{\vol (M)}{4n} \right) \geq \exp\left( \int_{M} \Omega(\varphi(p)) \, d{\vol}(p)\right) \quad \text{for
			all } \varphi \in R(M, n).
		\]
		It is clear that $R(M, n)\subset\mathrm{Map}_{\mu_n}\left( M, X, \rho: U_n, F_n,\frac{1}{4n} \right)$ by $e\in U_n$ and $\d_n$ small enough. Since every element of $S_n$ is $(\rho^{M}, \varepsilon)$-separated, one has
		
		\begin{align*}
			&\frac
			{1}{{\vol(M)}} \log \sum_{\varphi \in S_n} \exp\left( \int_{M} \Omega(\varphi(p)) \, d{\vol}(p) \right) \\
			&\leq\frac
			{1}{{\vol(M)}} \log	\left(|B|\sum_{\varphi \in R(M,n)} \exp\left( \int_M \Omega(\varphi(p)) \, d{\mathrm{vol}}(p) \right)\right)\\
			&\leq
			\frac
			{1}{{\vol(M)}} \log(|B|\cdot|R(M,n)|) + \int_{M} \Omega \, d\mu_n + \frac{1}{4n}
			\\
			&\leq \frac
			{1}{{\vol(M)}} \log\left( |B|\Sep_{\varepsilon}\left( \mathrm{Map}_{\mu_n}\left( M, X, \rho: U_n, F_n,\frac{1}{4n} \right) \right) \right) + \int \Omega \, d\mu_n + \frac{1}{4n}
			.
		\end{align*}
		
		Thus
		\begin{align*}
			&\frac{1}{\vol(M)} \log	M_{\Sigma}^{\ep}(\Omega, M, X, \rho: U_n, \delta_n)\\
			&\leq\frac{1}{\vol(M)} + \frac{1}{\vol(M)}\log\left( \sum_{\varphi \in S_n} \exp\left( \int_{M} \Omega(\varphi(p)) \, d{\vol}(p) \right) \right) \\
			&\leq \frac{1}{\vol(M)} + \frac{1}{\vol(M)}\log\left( |B|\Sep_{\varepsilon}\left( \mathrm{Map}_{\mu_n}\left( M, X, \rho: U_n, F_n,\frac{1}{4n} \right) \right) \right)+ \int \Omega \, d\mu_n + \frac{1}{4n}
			.
		\end{align*}

		Therefore, there exist $\mu_n \in B$ and a sequence $i_1 < i_2 < \cdots$ in $\mathbb{N}$ with
		\[
		P_{\Sigma}^{\ep}(\Omega, X, G, \rho, U_n, \delta_n) = \lim_{k \to \infty} \frac{1}{\vol(M_{i_k})} \log 	M_{\Sigma}^{\ep}(\Omega, M_{i_k}, X, \rho: U_n, \delta_n).
		\]
		such that
		\begin{align*}
			&\frac{1}{\vol(M_{i_k})} \log	M_{\Sigma}^{\ep}(\Omega, M_{i_k}, X, \rho: U_n, \delta_n)\\ &\leq \frac{1}{\vol(M_{i_k})}\log\left( |B|\Sep_{\varepsilon}\left( \mathrm{Map}_{\mu_n}\left( M_{i_k}, X, \rho: U_n,F_n, \frac{1}{4n} \right) \right) \right) +\frac{1}{\vol(M_{i_k})} + \int \Omega \, d\mu_n + \frac{1}{4n},
		\end{align*}
		
		for all $k \in \mathbb{N}$ and $|\int f\circ \gamma \ d\mu_n  - \int f \ d\mu_n | < 1/n$ for any $\gamma \in U_n$, and $f \in F_n$. Then

		\begin{equation}
			\begin{aligned}
				&P_{\Sigma}^{\varepsilon}(\Omega, X, G, \rho)\\
				&\leq P_{\Sigma}^{\varepsilon}(\Omega, X, G, \rho, U_n, \delta_n) \\
				&= \lim_{k \to \infty} \frac{1}{\vol(M_{i_k})} \log M_{\Sigma}^{\varepsilon}(\Omega, M_{i_k}, X, \rho: U_n, \delta_n) \\
				&\leq \lim_{k \to \infty} \left( \frac{1}{\vol(M_{i_k})} \log \left( |B_n| \Sep_{\varepsilon} \left( \mathrm{Map}_{\mu_n} \left( M_{i_k}, X, \rho: U_n, F_n, \frac{1}{4n} \right) \right) \right) + \frac{1}{\vol(M_{i_k})} \right) \\
				&\quad + \int \Omega \, d\mu_n + \frac{1}{4n} \\
				&\leq h_{\Sigma,\mu_n}^{\varepsilon} \left( X, \rho: U_n, F_n, \frac{1}{4n} \right) + \int \Omega \, d\mu_n + \frac{1}{4n}.
			\end{aligned}
			\label{eq:pressure-inequality}
		\end{equation}
		
		%
		%
		%
		%
		%
		%
		%
		Note that the sequence $\{\mu_n\}_{n=1}^\infty$ depends on $\varepsilon$. Let $\mu^\varepsilon$ be a weak* limit point of the sequence $\{\mu_n\}_{n=1}^\infty$ and we will show that $\mu^\varepsilon$ is a $G$-invariant measure.

		Let $\varphi \in \mathrm{Map}(M, X, \rho: U_n, \delta_n)$, $f \in F_n$ and $g \in U_n$.	We claim that	if $M$ is a sufficiently good approximation,  then
		\[
	\left| \int \left( f \circ \varphi(p) - f \circ g \circ \varphi(p) \right) d\ov{\vol}(p) \right| \leq \frac{7}{n}.
	\]
		\begin{proof}[Proof of claim]
			
			Denote $t = \left| \int_M \left( f(\varphi(p)) - f(g\varphi(p)) \right) \, \mathrm{dvol}(p) \right|$. To establish this inequality, we decompose $\alpha$ into two parts:
			\begin{align*}
				t_1 &= \left| \int_M f(\varphi(p)) \, d\ov{\vol}(p) - \int_{M[g]} f(\varphi(g. p)) \, d\ov{\vol}(p) \right|, \\
				t_2 &= \left| \int_{M[g]} f(\varphi(g.p)) \, d\ov{\vol}(p) - \int_M f(g\varphi(p)) \, d\ov{\vol}(p) \right|,
			\end{align*}
			where $M[g]$ denotes the subset of $M$ consisting of points $p$ for which the partial left-action $g.p$ is well-defined.
			
			We begin by estimating $t_1$. Applying Lemma \ref{lmp}, we obtain the identity:
			\[
			\int_{M[U_n]} f(\varphi(p)) \, d\ov{\vol}(p) = \int_{g^{-1}.M[U_n]} f(\varphi(g.p)) \, d\ov{\vol}(p).
			\]
			Since $M$ is an $(U_n, \delta_n)$-sofic approximation, we have $\mathrm{vol}(M[U_n]) \geq 1 - \delta_n$. Consequently,
			\[
			t_1= \left| \int_{M \setminus M[U_n]} f(\varphi(p)) \, d\ov{\vol}(p) - \int_{M[g] \setminus g^{-1} \cdot M[U_n]} f(\varphi(g.p)) \, d\ov{\vol}(p) \right| \leq 2\delta_n \left( \sup_{f \in F_n} \| f \|_\infty \right).
			\]
			
			Next, we estimate $	t_2$. Let $M'$ the set of all $p \in M[U_n]$ satisfying $\rho(\varphi(g.p), g\varphi(p)) < \sqrt{\d_n}$. By our choice of $\delta_n$, we have $|f(g\varphi(p)) - f(\varphi(g.p))| < \frac{1}{n}$ for all $p \in M'$. Therefore,
			\[
			\left| \int_{M'} \left( f(\varphi(g.p)) - f(g\varphi(p)) \right) \, d\ov{\vol}(p) \right| < \frac{1}{n}.
			\]
			
			Applying Markov's inequality and noting that $M' \subset M[U_n]$, we find $\ov{\mathrm{vol}}(M') \geq 1 - \delta_n - \sqrt{\delta_n}$. Thus,
		\begin{align*}
				t_2& \leq \frac{1}{n} + \left| \int_{M[g] \setminus M'} f(\varphi(g.p)) \, d\ov{\vol}(p) - \int_{M \setminus M'} f(g\varphi(p)) \, d\ov{\vol}(p) \right|\\
				& \leq \frac{1}{n} + 2(\delta_n + \sqrt{\delta_n}) \left( \sup_{f \in F_n} \| f \|_\infty \right).
		\end{align*}

			Combining these estimates, we obtain:
			\[
			t \leq t_1 + t_2 \leq \frac{1}{n} + (4\delta_n + 2\sqrt{\delta_n}) \left( \sup_{f \in F_n} \| f \|_\infty \right) < \frac{7}{n},
			\]
			which completes the proof.
		\end{proof}

		Hence, for any $g \in U_n$, $f \in F_n$, and every $\mu_\varphi$, using \eqref{eq:4.1} we conclude that
		\begin{equation} \label{eq:4.3}
			\begin{split}
				&\left| \int (f \circ g - f)  d\mu_\varphi \right| \\
				&\quad \leq \left| \int (f \circ g - f \circ {a_g})  d\mu_\varphi \right| + \left| \int f \circ {a_g} \, d\mu_\varphi - \int f \circ {a_g} \circ \varphi \, d\ov{\vol} \right| \\
				&\quad \quad + \left| \int (f \circ {a_g} \circ \varphi - f \circ \varphi)  d\ov{\vol} \right| + \left| \int f \circ \varphi  \,d\ov{\vol} - \int f  d\mu_\varphi \right| \\
				&\quad < \frac{1}{4n} + \frac{1}{4n} + \frac{7}{n} + \frac{1}{4n} = \frac{31}{4n}.
			\end{split}
		\end{equation}
		
		For any $f \in \bigcup_{n=1}^\infty F_n$ and $g \in G$,  we have
		\begin{align*}
			&\left| \int (f \circ g - f)  d\mu^\varepsilon \right|\\
			& \leq \left| \int f \circ g \, d\mu^\varepsilon - \int f \circ g\,  d\mu_n \right| + \left| \int (f \circ g - f)  d\mu_n \right| + \left| \int f  d\mu_n - \int f  d\mu^\varepsilon \right|.
		\end{align*}
		Then by \eqref{eq:4.3} and  the definition of weak* limits, we obtain
		\[
		\int (f \circ g - f)  d\mu^\varepsilon = 0.
		\]
		Furthermore, by our choice, $\bigcup_{n=1}^\infty F_n$ is dense in $C(X)$, and thus $\mu^\varepsilon$ is $G$-invariant.
		
		Let $F$ be a finite subset of $C(X)$, $U\subset G$ a precompact neighborhood of the identity, and $\varepsilon, \delta > 0$ are given. Take a large number $n_0$ with $U \subset U_{n_0}$, $1/n_0 \leq \delta$, and such that every $f \in F$ can be approximated by $f' \in F_{n_0}$ with $\| f' - f \|_\infty \leq \delta/4$, and moreover
		\[
		\sup_{f' \in F_{n_0}} \left| \int f'  d\mu_{n_0} - \int f'  d\mu^\varepsilon \right| \leq \delta/4.
		\]
		Then for any $\varphi \in \Map_{\mu_{n_0}}(M,X,\rho: U_{n_0}, F_{n_0}, \frac{1}{4n_0})$, we have
		\begin{align*}
			&\left| \int f \circ \varphi \, d\ov{\vol} - \int f  d\mu^\varepsilon \right| \\
			&\quad \leq \left| \int (f \circ \varphi - f' \circ \varphi)  d\ov{\vol} \right| + \left| \int f' \circ \varphi  d\ov{\vol} - \int f'  d\mu^\varepsilon \right| + \left| \int f'  d\mu^\varepsilon - \int f  d\mu^\varepsilon \right| \\
			&\quad \leq \| f' - f \|_\infty  + \left| \int f' \circ \varphi \, d\ov{\vol} - \int f'  d\mu_{n_0} \right| + \left| \int f'  d\mu_{n_0} - \int f'  d\mu^\varepsilon \right| +\| f' - f \|_\infty \\
			&\quad \leq \frac{\delta}{4} + \frac{1}{4n_0} + \frac{\delta}{4}+\frac{\delta}{4} \leq \delta.
		\end{align*}
		Thus, we have the inclusion $$\Map_{\mu_{n_0}}(M,X,\rho: U_{n_0}, F_{n_0}, \frac{1}{4n_0})\subset\Map_{\mu^\ep}(M,X,\rho: U, F, \d).$$
		Therefore using \eqref{eq:pressure-inequality} we have
		\begin{align*}				
			h_{\Sigma,\mu^\ep}^{\ep}\left( X, \rho: U,F, \delta \right) + \int\Omega \, d\mu_n&\ge		
			h_{\Sigma,\mu_{n_0}}^{\ep}\left( X, \rho: U_{n_0},F_{n_0}, \frac{1}{4{n_0}} \right) + \int\Omega \, d\mu_{n_0}-\frac{\delta}{4}	\\
			&\ge P_{\Sigma}^{\ep}(\Omega, X, G, \rho)-\frac{1}{4n_0}-\frac{\delta}{4}	\\
			&\ge P_{\Sigma}^{\ep}(\Omega, X, G, \rho)-\frac{\delta}{2}.	
		\end{align*}
		Since $U$, $F$, $\delta$, $\varepsilon$ are arbitrary, we conclude that
		\[
		P_{\Sigma}^{\ep}(\Omega, X, G, \rho)\leq h_{\Sigma,\mu^\ep}^{\ep}(\rho) + \int\Omega \, d\mu_n\leq h_{\Sigma,\mu^\ep}(\rho) + \int\Omega \, d\mu_n
		\]
		and therefore
		
		\begin{align*}	
			P_{\Sigma}(\Omega, X, G)= \sup_{\ep>0} P_{\Sigma}^{\ep}(\Omega, X, G, \rho)&\leq\sup_{\ep>0} \left(h_{\Sigma,\mu^\ep}(\rho) + \int\Omega \, d\mu^\ep\right)\\
			&\leq\sup_{\mu\in\pr(X,G)} \left(h_{\Sigma,\mu}(\rho) + \int\Omega \, d\mu\right).
		\end{align*}
		
	\end{proof}
	
	Now we can prove the variational principle.
	
	\begin{proof}[Proof of Theorem~\ref{main3}]
		Let $\rho$ be a compatible metric on $X$, $U\subset G$ a precompact neighborhood of the identity. Let $\mu \in \pr(X,G)$ and $\delta, \varepsilon > 0$. Put $F_1 = \{\Omega\}$. Fix $i \in \mathbb{N}$. Let $S_i$ be a $(\rho^{M_i}, \varepsilon)$-separated subset of $\mathrm{Map}_{\mu}(M_{i}, X, \rho: U,F_1, \delta)$ with maximal cardinality. Then $S_i$ is also a $(\rho^{M_i}, \varepsilon)$-separated subset of $\mathrm{Map}(M_{i}, X, \rho: U, \delta)$.
		
		Since the function $x \mapsto \log x$ for $x > 0$ is concave, one has
		\[
		\log \sum_{\varphi \in S_i} \frac{1}{|S_i|}\exp\left( \int_{M_i} \Omega(\varphi(p)) \, d{\vol}(p) \right) \geq \frac{1}{|S_i|} \sum_{\varphi \in S_i} \int_{M_i} \Omega(\varphi(p)) \, d{\vol}(p).
		\]
		Hence
		\begin{align*}	
			\log \sum_{\varphi \in S_i} \exp\left( \int_{M_i} \Omega(\varphi(p)) \, d{\vol}(p) \right)
			&\geq \log |S_i| + \frac{1}{|S_i|} \sum_{\varphi \in S_i} \int_{M_i} \Omega(\varphi(p)) \, d{\vol}(p)\\
			&\geq \log |S_i| + \frac{1}{|S_i|} \sum_{\varphi \in S_i} \left( \int_X \Omega \, d\mu - \delta \right)\vol(M_i)\\
			&= \log |S_i| + \left( \int_X \Omega \, d\mu - \delta \right) \vol(M_i).
		\end{align*}
		Thus for all  precompact neighborhood of the identity $U\subset G$ and all $\delta, \varepsilon > 0$,
		\[
		P_{\Sigma}^\ep(\Omega, X, G, \rho, U, \delta) + \delta \geq \limsup_{i \to \infty} \frac{1}{\vol(M_i)} \log |S_i| + \int_X \Omega \, d\mu.
		\]
		Then
		\[
		P_{\Sigma}^\ep(\Omega, X, G, \rho, U) \geq \inf_{\delta > 0}\limsup_{i \to \infty} \frac{1}{\vol(M_i)} \log |S_i| + \int_X \Omega \, d\mu.
		\]
		Since $U$ and $\varepsilon$ are arbitrary
		\[
		P_\Sigma(\Omega, X, G) \geq h_{\Sigma,\mu}(X) + \int_X \Omega \, d\mu.
		\]
		Combining with Lemma~\ref{lma3}, we get
		\[
		P_\Sigma(\Omega, X, G) = \sup \left\{ h_{\Sigma,\mu}(X) + \int_X \Omega \, d\mu : \mu \in \pr(X,G) \right\}.
		\]
	\end{proof}

	\subsection{Properties of sofic topological pressure}
	\begin{prop}\label{prop}
		If $\Omega, \Lambda \in C(X)$, $g \in G$ and $c \in \mathbb{R}$ then the following are true.
		\begin{enumerate}
			\item[(i)] $P_\Sigma(0, X, G) = h_\Sigma(X)$.
			\item[(ii)] $P_\Sigma(\Omega + c, X, G) = P_\Sigma(\Omega, X, G) + c$.
			\item[(iii)] $P_\Sigma(\Omega + \Lambda, X, G) \leq P_\Sigma(\Omega, X, G) + P_\Sigma(\Lambda, X, G)$.
			\item[(iv)] $\Omega \leq \Lambda$ implies $P_\Sigma(\Omega, X, G) \leq P_\Sigma(\Lambda, X, G)$. In particular, 
			$h_\Sigma(X) + \min \Omega \leq P_\Sigma(\Omega, X, G) \leq h_\Sigma(X) + \max \Omega$.
			\item[(v)] $P_\Sigma(\cdot, X, G)$ is either finite valued or constantly $\pm\infty$.
			\item[(vi)] If $P_\Sigma(\cdot, X, G) \neq \pm\infty$, then 
			$|P_\Sigma(\Omega, X, G) - P_\Sigma(\Lambda, X, G)| \leq \|\Omega - \Lambda\|_\infty$,
			where $\|\cdot\|_\infty$ is the supremum norm on $C(X)$.
			\item[(vii)] If $P_\Sigma(\cdot, X, G) \neq \pm\infty$ then $P_\Sigma(\cdot, X, G)$ is convex.
			\item[(viii)] $P_\Sigma(\Omega + \Lambda \circ g - \Lambda, X, G) = P_\Sigma(\Omega, X, G)$.
			\item[(ix)] $P_\Sigma(c\Omega, X, G) \leq c \cdot P_\Sigma(\Omega, X, G)$ if $c \geq 1$ and 
			$P_\Sigma(c\Omega, X, G) \geq c \cdot P_\Sigma(\Omega, X, G)$ if $c \leq 1$.
			\item[(x)] $|P_\Sigma(\Omega, X, G)| \leq P_\Sigma(|\Omega|, X, G)$.
		\end{enumerate}
	\end{prop}
	
	\begin{proof}
		Let $\rho$ be a compatible metric on $X$. Let $\varepsilon, \delta > 0$ and $U\subset G$ a precompact neighborhood of the identity.
		
		\begin{enumerate}
			\item[(i)--(iv)] These are clear from the definition of sofic topological pressure and Remark \ref{rem1}.
			
			\item[(v)] From (i) and (ii) we get $P_\Sigma(\Omega, X, G) = \pm\infty$ if and only if $h_\Sigma(X) = \pm\infty$.
			
			\item[(vi)] Follows from (iii) and (iv).
			
			\item[(vii)] By Hölder's inequality, if $p \in [0, 1]$ and $E$ is a finite subset of $\mathrm{Map}(M, X, \rho: U, \delta)$ then we have
			\begin{align*}
				\sum_{\varphi \in E} \exp\left[p \int_{M} \Omega(\varphi(p)) \, d{\vol}(p) + (1-p) \int_{M} \Lambda(\varphi(p)) \, d{\vol}(p)\right] \\
				\leq \left[\sum_{\varphi \in E} \exp\left(\int_{M} \Omega(\varphi(p)) \, d{\vol}(p)\right)\right]^p 
				\left[\sum_{\varphi \in E} \exp\left(\int_{M} \Lambda(\varphi(p)) \, d{\vol}(p)\right)\right]^{1-p}.
			\end{align*}
			Therefore,
			\begin{align*}
				M^{\varepsilon}_{\Sigma}(p\Omega + (1-p)\Lambda, X, G, \rho, U, \delta) \\
				\leq M^{\varepsilon}_{\Sigma}(\Omega, X, G, \rho, U, \delta)^p 
				\cdot M^{\varepsilon}_{\Sigma}(\Lambda, X, G, \rho, U, \delta)^{1-p},
			\end{align*}
			and (vii) follows.
			
			\item[(viii)] It is clear from the variational principle.

			\item[(ix)] If $a_1, \dots, a_k$ are positive numbers with $\sum_{i=1}^k a_i = 1$ then $\sum_{i=1}^k a_i^c \leq 1$ when $c \geq 1$, and $\sum_{i=1}^k a_i^c \geq 1$ when $c \leq 1$. Hence if $b_1, \dots, b_k$ are positive numbers then
			\begin{align*}
				\sum_{i=1}^k b_i^c \leq \left(\sum_{i=1}^k b_i\right)^c \quad \text{when } c \geq 1, \\
				\sum_{i=1}^k b_i^c \geq \left(\sum_{i=1}^k b_i\right)^c \quad \text{when } c \leq 1.
			\end{align*}
			Therefore, if $E$ is a finite subset of $\mathrm{Map}(M, X, \rho: U, \delta)$ we have
			\begin{align*}
				\sum_{\varphi \in E} \exp\left(c \int_{M} \Omega(\varphi(p)) \, d{\vol}(p) \right)
				\leq \left(\sum_{\varphi \in E} \exp\left(\int_{M} \Omega(\varphi(p)) \, d{\vol}(p)\right)\right)^c \quad \text{when } c \geq 1, \\
				\sum_{\varphi \in E} \exp\left(c \int_{M} \Omega(\varphi(p)) \, d{\vol}(p) \right)
				\geq \left(\sum_{\varphi \in E} \exp\left(\int_{M} \Omega(\varphi(p)) \, d{\vol}(p)\right)\right)^c \quad \text{when } c \leq 1.
			\end{align*}
			Then (ix) follows.
			
			\item[(x)] From (iv) and (ix) we get (x). \qedhere
		\end{enumerate}
	\end{proof}
	
	Let $\mathcal{B}_X$ be the $\sigma$-algebra of Borel subsets of $X$. Recall that a finite signed measure is a map $\mu: \mathcal{B}_X \to \mathbb{R}$ satisfying
	\[
	\mu\left(\bigcup_{i=1}^\infty A_i\right) = \sum_{i=1}^\infty \mu(A_i),
	\]
	whenever $\{A_i\}_{i=1}^\infty$ is a pairwise disjoint collection of members of $\mathcal{B}_X$.
	
	Using  the properties of sofic topological pressure, we can also give a sufficient condition for a finite signed measure to be a member of $\pr(X,G)$. See~{\cite[Theorem 9.11]{W75} for $\Z$-actions, {\cite[Theorem 6.2]{C13} for countable discrete sofic group actions and the proofs holds for locally compact sofic group actions naturally. For the sake of completeness, we provide the proofs.
			\begin{thm}\label{signed}
				Let $\Sigma = (M_i)_{i=1}^{\infty}$ be a sofic approximation to $G$, $(X)$ be a $G$-system and $h_\Sigma(X) \neq \pm \infty$. Let $\mu: \mathcal{B}_X \to \mathbb{R}$ be a finite signed measure. If $\int_X \Omega \, d\mu \leq P_\Sigma(\Omega, X, G)$ for all $\Omega\in C(X)$, then $\mu $ is a $G$-invariant measure.
			\end{thm}
			
			\begin{proof}
				We prove the theorem in three steps.
				
				\textbf{Step 1: Prove that $\mu$ is non-negative.} 
				Let $\Omega  \in C(X)$ with $\Omega \geq 0$. We aim to show that $\int_X \Omega  \, d\mu \geq 0$. For any $\ep> 0$ and integer $n > 0$, consider the function $n(\Omega  + \ep)$. By the assumption on $\mu$, we have
				\[
				\int_X n(\Omega +\ep) \, d\mu = - \int_X -n(\Omega +\ep) \, d\mu \geq -P_\Sigma(-n(\Omega +\ep), X, G).
				\]
				By Proposition \ref{prop} (iv), we have the bound
				\[
				P_\Sigma(-n(\Omega +\ep), X, G) \leq h_\Sigma(X) + \max(-n(\Omega +\ep)) = h_\Sigma(X) - n \min(\Omega +\ep).
				\]
				Since $\Omega \geq 0$ and $\ep > 0$, we have $\min(\Omega +\ep) \geq \ep > 0$. Thus,
				\[
				- P_\Sigma(-n(\Omega +\ep), X, G) \geq -h_\Sigma(X) + n \min(\Omega +\ep) \geq -h_\Sigma(X) + n \ep.
				\]
				Combining the inequalities, we obtain
				\[
				\int_X n(\Omega +\ep) \, d\mu \geq -h_\Sigma(X) + n\ep.
				\]
				For sufficiently large $n$, the right-hand side is positive, so $\int_X n(\Omega +\ep) \, d\mu > 0$. Dividing by $n > 0$, we get $\int_X (\Omega +\ep) \, d\mu > 0$. Since $\ep > 0$ is arbitrary, it follows that $\int_X \Omega \, d\mu \geq 0$. Therefore, $\mu$ is a non-negative measure.
				
				\textbf{Step 2: Prove that $\mu(X) = 1$.}
				Consider the constant function $\Omega = n$ for $n \in \mathbb{Z}$. By assumption,
				\[
				\int_X n \, d\mu = n \mu(X) \leq P_\Sigma(n, X, G).
				\]
				By Proposition \ref{prop} (iv), we have $P_\Sigma(n, X, G) = h_\Sigma(X) + n$. Thus,
				\[
				n \mu(X) \leq h_\Sigma(X) + n.
				\]
				We now consider two cases:
				\begin{itemize}
					\item If $n > 0$, dividing both sides by $n$ gives $\mu(X) \leq 1 + \frac{h_\Sigma(X)}{n}$. Taking the limit as $n \to \infty$, we obtain $\mu(X) \leq 1$.
					\item If $n < 0$, dividing both sides by $n$ (which is negative) $\mu(X) \geq 1 + \frac{h_\Sigma(X)}{n}$. Taking the limit as $n \to -\infty$, we obtain $\mu(X) \geq 1$.
				\end{itemize}
				Therefore, $\mu(X) = 1$.
				
				\textbf{Step 3: Prove that $\mu$ is $G$-invariant.}
				Let $g \in G$ and $f \in C(X)$. We need to show that $\int_X \Omega \circ g \, d\mu = \int_X \Omega \, d\mu$. For any integer $n$, by assumption,
				\[
				\int_X n(\Omega \circ g - \Omega) \, d\mu \leq P_\Sigma(n(\Omega \circ g - \Omega), X, G).
				\]
				By Proposition \ref{prop} (viii), we have $P_\Sigma(n(\Omega \circ g - \Omega), X, G) = h_\Sigma(X)$. Thus,
				\[
				n \int_X (\Omega \circ g - \Omega) \, d\mu \leq h_\Sigma(X).
				\]
				We analyze this inequality for different values of $n$:
				\begin{itemize}
					\item If $n > 0$, dividing both sides by $n$ gives $\int_X (\Omega \circ g - \Omega) \, d\mu \leq \frac{h_\Sigma(X)}{n}$. Taking the limit as $n \to \infty$, we get $\int_X (\Omega \circ g - \Omega) \, d\mu \leq 0$.
					\item If $n < 0$, dividing both sides by $n$ (which is negative)  $\int_X (\Omega \circ g - \Omega) \, d\mu \geq \frac{h_\Sigma(X)}{n}$. Taking the limit as $n \to -\infty$, we get $\int_X (\Omega \circ g - \Omega) \, d\mu \geq 0$.
				\end{itemize}
				Hence, $\int_X (\Omega \circ g - \Omega) \, d\mu = 0$, which implies $\int_X \Omega \circ g - \Omega \, d\mu = \int_X f \, d\mu$ for all $f \in C(X)$ and $g \in G$. Therefore, $\mu$ is $G$-invariant.
			\end{proof}
			
			\begin{ques}
				It is known that if $G$ is amenable and unimodular then it is soﬁc(see \cite{BB22}) and if  $G$ is countable amenable then the soﬁc topological pressure agrees with classical topological pressure(see \cite{C13}). Is it true that the result holds at the locally compact level?
			\end{ques}

			\bibliographystyle{amsalpha}

		\end{document}